# LARGE DEVIATION FOR DIFFUSIONS AND HAMILTON–JACOBI EQUATION IN HILBERT SPACES

By Jin Feng

*University of Massachusetts–Amherst*

Large deviation for Markov processes can be studied by Hamilton–Jacobi equation techniques. The method of proof involves three steps: First, we apply a nonlinear transform to generators of the Markov processes, and verify that limit of the transformed generators exists. Such limit induces a Hamilton–Jacobi equation. Second, we show that a strong form of uniqueness (the comparison principle) holds for the limit equation. Finally, we verify an exponential compact containment estimate. The large deviation principle then follows from the above three verifications.

This paper illustrates such a method applied to a class of Hilbert-space-valued small diffusion processes. The examples include stochastically perturbed Allen–Cahn, Cahn–Hilliard PDEs and a one-dimensional quasilinear PDE with a viscosity term. We prove the comparison principle using a variant of the Tataru method. We also discuss different notions of viscosity solution in infinite dimensions in such context.

**1. Introduction.** We are interested in large deviation for small randomly perturbed diffusion processes in a Hilbert state space $E$. When $E = R^d$, this is known as the Freidlin and Wentzell theory [23]. The proofs in [23] rely upon the Girsanov transformations. The idea is to estimate probability of an atypical, large deviant event under the given probability law through a change of measure, so that the event becomes most probable under the new law. Such technique is also repeatedly used in the Donsker and Varadhan theory [13] regarding occupation measures, which is another kind of large deviation concerning ergodic phenomena instead of small random perturbations.

There exists a different approach to the above mentioned large deviation problems. In the late 1970s, Fleming [19] introduced a logarithmic trans-form to generators of Markov processes, giving exit probabilities an optimal control interpretation. This observation allowed us to characterize the large deviation convergence for exit probabilities as convergence of solutions for a sequence of Hamilton–Jacobi equations. Later, Evans and Ishii [16], Fleming and Souganidis [21, 22], among others, applied the theory of viscosity solution to this context, enabling the approach to cover a wider variety of examples. In particular, this includes $R^d$-valued diffusions with vanishing stochastic terms. During early developments of this approach, the applicable settings and structural conditions required were relatively restrictive as compared to the Girsanov transformation approach. However, this can be fixed by refining techniques on the viscosity solution techniques and on the large deviation theory. Feng and Kurtz [18] recently carried out such a program which expands the theory.

The general setting in [18] allows the state space $E$ to be a metric space. One of the key technical conditions assumed is a strong form of uniqueness (i.e., the comparison principle, Definition 1.15) for a limit Hamilton–Jacobi type equation. For small perturbation type large deviations, such equation is usually a first-order nonlinear partial differential equation. When $E = R^d$, or a subset of it, the comparison principle can usually be verified by well-known criteria in PDE theory. Using these techniques, [18] treats both the classical Freidlin–Wentzell theory and the Donsker–Varadhan theory within one framework using the generator convergence approach.

When we study large deviation for stochastic PDEs or interacting particles, we usually encounter function- or measure-valued state space. Comparison principles of these types, however, are much less well understood. On the one hand, there exists an extensive PDE literature regarding first-order Hamilton–Jacobi equations in Hilbert/Banach spaces (e.g., [7, 8, 31, 32]). On the other hand, the operators derived in the large deviation context frequently exhibit subtle differences relative to those studied in the PDE literature. Indeed, in the case of applications to interacting particle systems, it is more natural to consider the state space as the space of probability measures, rather than a Banach space.

We restrict attention to Hilbert-space-valued diffusions only in this paper.

1.1. *Background.* We consider the large deviation for Hilbert-space-valued diffusions with a possibly nonlinear drift term. To outline the approach and identify difficulties ahead of us, first, we review a general result (adapted to the situation of this paper) developed in [18].

Let $X_n$, $n = 1, 2, \dots$, be a sequence of metric-space $S$-valued random variables. Varadhan and Brycs (e.g. Theorems 4.3.1 and 4.4.2 in [12]) discovered the following moment characterization of large deviation convergence.

PROPOSITION 1.1.

(a) *Suppose* $\{X_n\}$ *satisfies the large deviation principle (Definition 1.17) with a good rate function* $I$. *Then for each* $f \in C_b(S)$ *(bounded continuous functions on* $S$*), if we define* $\Lambda_n(f) = n^{-1} \log E[\exp\{nf(X_n)\}]$,

$$\lim_{n\to+\infty} \Lambda_n(f) = \lim_{n\to+\infty} \frac{1}{n} \log E[e^{nf(X_n)}] = \sup_{x\in S}\{f(x) - I(x)\} = \Lambda(f). \tag{1.1}$$

(b) *Suppose that* $\{X_n\}$ *is exponentially tight (Definition 1.17) and that the limit (1.1) exists for each* $f \in C_b(S)$. *Then* $\{X_n\}$ *satisfies the large deviation with good rate function*

$$I(x) = \sup_{f\in C_b(S)} (f(x) - \Lambda(f)). \tag{1.2}$$

See Theorems 4.3.1 and 4.4.2 in [12].

The main result in [18] can be viewed as a process version of the above theorem, expressed at an infinitesimal level.

To explain the result, we proceed informally first. Let $\{X_n(t), 0 \le t < +\infty; n = 1, 2, \dots\}$ denote a sequence of metric-space $E$-valued Markov processes. For simplicity, we assume the trajectories are continuous. By the Markov property and continuity of the trajectories, we expect the large deviation of $\{X_n\}$ follows from that of the transition probability measures $\{P(X_n(t) \in dy|X_n(0) = x)\}$ for all $x \in E$ and $t \ge 0$. By Proposition 1.1, we also expect this to be implied by convergence of the functionals

$$V_n(t)f(x) \equiv \frac{1}{n} \log E[e^{nf(X_n(t))}|X_n(0) = x] \to V(t)f(x) \qquad \text{for some } V(t)f,$$

where $f \in D \subset C_b(E)$, and $D$ is sufficiently dense in $C_b(E)$ in appropriate sense.

It turns out that, by the Markov property, $V_n$ forms a nonlinear operator semigroup

$$V_n(s)V_n(t) = V_n(t+s), \qquad s, t \ge 0.$$

Hence $\{V(t) : t \ge 0\}$, viewed as a collection of operators acting on functions, should form a semigroup as well. We identify the generator of $V_n$ next:

$$\begin{aligned} H_n f(x) &= \lim_{t\to 0+} \frac{1}{t}(V_n(t)f(x) - V_n(0)f(x)) \\ &= \lim_{t\to 0+} \frac{1}{t}\frac{1}{n} \log E[e^{n(f(X_n(t))-f(x))}|X_n(0) = x] \\ &= \frac{1}{n} e^{-nf(x)} A_n e^{nf}(x), \end{aligned}$$

where $A_n$ is the generator for the process $X_n$

$$A_n g(x) = \lim_{t\to 0+} \frac{1}{t} \log E[g(X_n(t)) - g(x)|X_n(0) = x].$$

The above transformation from $A_n$ to $H_n$ is essentially the logarithmic transform by Fleming [19]. If we denote by $H$ the generator for $V$, then we expect generator convergence $H_n \to H$ will imply semigroup convergence $V_n \to V$, which is suggested by the semigroup generation theorem:

$$V_n(t)h = \lim_{k\to+\infty}\left(I - \frac{t}{k}H_n\right)^{-k} h$$

and

$$V(t)h = \lim_{k\to+\infty}\left(I - \frac{t}{k}H\right)^{-k} h. \tag{1.3}$$

Going backward in the reasoning, modulo regularity conditions, we expect convergence $H_n \to H$ will give the large deviation of $\{X_n\}$.

There are practical problems if we want to rigorously apply the above program to examples. First, the formula in (1.3) requires

$$(I - \alpha H)f = h \tag{1.4}$$

to hold in the classical sense for all $h \in D$ and $\alpha > 0$ (if $H$ is dissipative, then such $f$ is also unique). However, this is extremely hard to verify for most examples. Therefore we are forced to modify the above formulation by using a type of weak solution called the *viscosity solution* (Definition 1.14). By weakening the type of solution needed for (1.4), we have to require a strong form of uniqueness condition known as the *comparison principle* (Definition 1.15). Informally, this principle states that, if upper semicontinuous $\overline{f}$ and lower semicontinuous $\underline{f}$ satisfy

$$(I - \alpha H)\overline{f} \le h \quad \text{and} \quad (I - \alpha H)\underline{f} \ge h,$$

then $\overline{f} \le \underline{f}$. The $\overline{f}$ and $\underline{f}$ are called, respectively, a subsolution and a supersolution. Existence of sub- and supersolutions, in the large deviation context here, can be constructed by generalizing a procedure due to Barles and Perthame [2, 3]. When $E$ is noncompact, in order for such argument to go through, we need another crucial condition (Condition 2.2) for the transition probabilities. Such condition requires the processes to be concentrated on a compact subset of the state space with high probability.

Note that the above formulation is only based on inequalities for sub- and supersolutions. This provides an opportunity for further relaxations on conditions. We can introduce two more operators: $H_0, H_1$ so that $Hf \le H_0 f$ and $Hf \ge H_1 f$ for all $f \in \mathcal{D}(H) \cap \mathcal{D}(H_i)$. Then

$$(I - \alpha H_0)\overline{f} \le h \quad \text{and} \quad (I - \alpha H_1)\underline{f} \ge h.$$

Suppose that the comparison principle still holds for the above two "inequations" (i.e., $\overline{f} \le \underline{f}$). The construction of $\underline{f}, \overline{f}$ by the Barles–Perthame procedure then reveals that $\overline{f} = \underline{f} = f \in C_b(E)$. Hence, each $h$ uniquely corresponds to an $f \in C_b(E)$, and we can denote such correspondence by $f = R_\alpha h$. Consequently, at least formally, $R_\alpha = (I - \alpha H)^{-1}$. In other words, $H_0$, $H_1$ implicitly determine $H$ through its resolvent, and $V(t)h = \lim_n R^n_{t/n} h \in C_b(E)$. We can now completely avoid using $H$ in the above program by replacing condition $H_n \to H$ by: for each $f \in \mathcal{D}(H_i)$,

$$H_1 f \le \liminf_n H_n f_n, \qquad \limsup_n H_n f_n \le H_0 f \qquad \text{some } f_n \to f.$$

The above generalization is useful for applications where $E$ is infinite dimensional. We illustrate this next.

In general, the comparison principle proof relies upon test functions which behave like distance functions. For instance, in the case $E = R^d$, these test functions take the form $f(x) = (\mu/2)|x - y|^2$, $\mu \in R$ (see [5]). For $H_n$'s which are differential operators, there is no difficulty to include such functions in the domain. Furthermore, identifying $Hf$ as a limit of $H_n f$ is usually straightforward.

However, the situation becomes tricky when $E$ is infinite dimensional (e.g., Examples 1.2, 1.5 and 1.8). For instance, let $E = L^2(\mathcal{O})$ and $\mathcal{O} = [0, 1)$ with periodic boundary, and

$$H_n f(x) = \langle \Delta x, Df(x)\rangle + \tfrac{1}{2}\|Df(x)\|^2 + o_f(1),$$

where $o_f$ may depend on $f$. See (1.33) for definition of $Df$. We expect

$$Hf(x) = \langle \Delta x, Df(x)\rangle + \tfrac{1}{2}\|Df(x)\|^2.$$

But then, even for

$$f(x) = (\mu/2)\|x - y\|^2, \qquad \mu \in R, \tag{1.5}$$

where $y$ is arbitrarily smooth,

$$\lim_{x_n \to x} H_n f(x_n) \ne Hf(x).$$

Note that in this case, $Df(x) = \mu(x - y)$ and $\langle \Delta x, Df(x)\rangle$ is well defined as a function taking value in extended reals

$$\langle \Delta x, Df(x)\rangle = \mu\langle \Delta x, x - y\rangle = -\mu\|\nabla x\|^2 + \mu\langle \nabla x, \nabla y\rangle.$$

Assuming $\mu > 0$, by lower semicontinuity of $\|\nabla x\|$, we can however obtain $\limsup_{x_n \to x} H_n f(x_n) \le Hf(x)$. Similarly, by reversing the inequality, we can verify a lower bound estimate for the case $\mu < 0$.

More generally, if the above $\Delta$ is replaced by a general nonlinear dissipative operator $C$ [assuming the domain of $C$ is not the entire $E$, $\mathcal{D}(C) \ne E$], integration by parts may not even make sense any more. Consequently

$$Hf(x) = \langle Cx, Df(x)\rangle + \tfrac{1}{2}\|Df(x)\|^2$$

does not make sense for all $x \in E$, even if $y \in \mathcal{D}(C)$. Note that, if $C$ is dissipative,

$$\mu\langle Cx, x-y\rangle \le \mu\langle Cy, x-y\rangle \qquad \forall\, x, y \in \mathcal{D}(C), \mu > 0.$$

The right-hand side of the above is continuous in $x$, and can be extended to all $x \in E$ easily. Hence at least for $f$ of the form (1.5) with $\mu > 0$, if $H_n f(x_n)$ is defined, then it can be estimated from above by

$$\begin{aligned}\limsup_{x_n \to x} H_n f(x_n) &\le \lim_n \mu\langle Cy, x_n - y\rangle + \tfrac{1}{2}\|Df(x_n)\|^2 + o_f(1)\\ &= \mu\langle Cy, x-y\rangle + \tfrac{1}{2}\|Df(x)\|^2 \equiv H_0 f(x).\end{aligned}$$

Furthermore, such $H_0 f \in C(E)$ and is everywhere well defined. By the arbitrariness of $y \in \mathcal{D}(C)$, we hope such $H_0$ provides a sharp estimate on the asymptotics of $H_n$'s. Similarly, we can also estimate $H_n f$ from below by some $H_1 f$, if $\mu < 0$.

We call (1.4) a Hamilton–Jacobi equation, because of its connection with the optimal control problem. By the Markov property on the processes $X_n$, the $A_n$'s satisfy the maximum principle. The $H_n$'s, obtained as a transform of the $A_n$'s, also satisfy a nonlinear maximum principle. So does the limiting $H$ (and frequently, $H_0, H_1$). Using this property, the following variational representation of $H_n$ can usually be proved [17]:

$$H_n f(x) = \sup_{u \in U}(B_n f(x,u) - L_n(x,u)),$$

where $U$ is some auxiliary metric space, and for each $u$ fixed, $B_n f(\cdot, u)$ is a linear operator satisfying the maximum principle in $C_b(E)$, $L_n$ is a lower semicontinuous bivariate function. In the limit, $H$ is supposed to have a similar structure. This is known as the Nisio representation of generator for Hamiltonian operator $H$ in optimal control theory [25]. Based upon such representation, [18] proved theorems ensuring a simpler, variational representation of the rate function in an "action integral" form.

1.2. *Basic setup.* Let Hilbert-space-valued diffusion processes

$$dX_n(t) = \hat{C}_n X_n(t)\,dt + \frac{1}{\sqrt{n}} B_n(X_n(t))\,dW(t), \tag{1.6}$$

where $W$ is a cylindrical Wiener process [see (1.28)] on a separable real Hilbert space $U_0$, $E$ is another separable real Hilbert space, $\hat{C}_n - \omega I$ is an $m$-dissipative (possibly) nonlinear operator on $E$ for some $\omega > 0$: that is,

$$\langle \hat{C}_n x - \hat{C}_n y, x - y\rangle \le \omega\|x-y\|^2 \qquad \forall\, x, y \in \mathcal{D}(\hat{C}_n),$$

and the range of $I - \alpha\hat{C}_n$ satisfies

$$\mathcal{R}(I - \alpha\hat{C}_n) = E \qquad \forall\, \alpha > 0.$$

We also assume that $0 \in \mathcal{D}(\hat{C}_n)$, and $B_n(x): U_0 \to E$ is a Hilbert–Schmidt operator for each $x \in E$ fixed. More conditions are needed in order to make sense of the solution and large deviation result of (1.6); we delay them until the statement of the respective theorem. To simplify the presentation, we will actually deal with another form of the above equation:

$$dX_n(t) = C_n X_n(t)\, dt + F_n(X_n(t))\, dt + \frac{1}{\sqrt{n}} B_n(X_n(t))\, dW(t), \tag{1.7}$$

where $C_n$ is $m$-dissipative, $C_n 0 = 0$ and $F_n(x): E \to E$ is globally Lipschitz in $x$. To rewrite (1.6) into the form of (1.7), we take

$$C_n x = \hat{C}_n x - \hat{C}_n 0 - \omega x, \qquad F_n(x) = \hat{C}_n 0 + \omega x.$$

We provide three examples.

Example 1.2 (Stochastic Allen–Cahn equation). Let $\mathcal{O} = [0,1)^d$, $d = 1, 2, 3, \ldots$, with periodic boundary condition; we associate $L^2(\mathcal{O})$ with the usual inner product

$$\langle x, y\rangle \equiv \int_{\theta \equiv (\theta_1, \ldots, \theta_d) \in \mathcal{O}} x(\theta) y(\theta)\, d\theta, \qquad x, y \in L^2(\mathcal{O}).$$

By a stochastic Allen–Cahn equation, we refer to the following formally written stochastic PDE

$$\begin{aligned} \frac{\partial}{\partial t} Y_n(t,\theta) &= \Delta Y_n(t,\theta) - V'(Y_n(t,\theta)) \\ &\quad + \frac{1}{\sqrt{n}} \sigma(\theta; Y_n(t)) \frac{\partial^{d+1}}{\partial t\, \partial\theta_1 \cdots \partial\theta_d} \beta(t,\theta), \end{aligned} \tag{1.8}$$

where $\beta(t,\theta)$ is a Brownian sheet over space–time $(t,\theta) \in [0,\infty) \times \mathcal{O}$, $V \in C^1(R)$ and

$$\sigma(\theta, y) = \varphi(\theta, \langle y, \xi_1\rangle, \ldots, \langle y, \xi_k\rangle) \tag{1.9}$$

for some $\xi_1, \ldots, \xi_k \in L^2(\mathcal{O})$ and $\varphi(\theta, r_1, \ldots, r_k): R^{k+1} \to R$.

The above is a stochastically perturbed reaction–diffusion type equation. Among other applications, it has been used in material science as a phenomenological model of material interface movements due to molecular-level adsorption–desorption processes. $V \in C^1(R)$ is usually a double- or multiple-well potential function. From large deviations for $\{Y_n\}$, we can extract information about metastability of the whole system when the temperature is small.

It is well known that, in dimension $d \geq 2$, (1.8) admits no $L^2(\mathcal{O})$-valued solution. We will actually consider an approximate version of it which is defined on truncated Fourier modes. The number of modes goes to infinity as $n \to +\infty$. Such consideration is motivated by the fact that, in the above mentioned application, (1.8) should only be viewed as a formal limit for some stochastic Ginzburg–Landau equation defined on finite lattices or on truncated Fourier modes [30]. It is usually the rescaling limits of these finite systems which we really care about, rather than the continuum level (1.8). We mention that large deviation for the lattice case is studied in [18].

To rigorously define the processes, we let

$$\phi_1(r) \equiv 1, \qquad \phi_{2k-1}(r) \equiv \sqrt{2}\cos(2\pi k r),$$
$$\phi_{2k}(r) \equiv \sqrt{2}\sin(2\pi k r), \qquad k = 1, 2, \ldots,$$

and

$$\mu_1 \equiv 0, \qquad \mu_{2k-1} = \mu_{2k} \equiv 4\pi^2 k^2.$$

It follows that

$$-\phi_j'' = \mu_j \phi_j. \tag{1.10}$$

Therefore

$$\begin{aligned} \{e_k(\theta) &\equiv \phi_{k_1}(\theta_1) \times \phi_{k_2}(\theta_2) \times \cdots \times \phi_{k_d}(\theta_d), \\ k &\equiv (k_1, \ldots, k_d), k_j = 1, 2, \ldots; j = 1, \ldots, d\} \end{aligned} \tag{1.11}$$

forms a complete orthonormal basis for $E \equiv U_0 \equiv L^2(\mathcal{O})$. Denote

$$\lambda_k \equiv \mu_{k_1} \times \cdots \times \mu_{k_d}; \tag{1.12}$$

then

$$-\Delta e_k = \lambda_k e_k.$$

Let $\{\beta_k(t), k \equiv (k_1, \ldots, k_d), k_j = 1, \ldots, j = 1, \ldots, d\}$ be a sequence of i.i.d. real-valued standard Brownian motion, and let

$$\beta(t, \theta) \equiv \sum_k \beta_k(t) \int_{r_1=0}^{\theta_1} \cdots \int_{r_d=0}^{\theta_d} e_k(r_1, \ldots, r_d)\, dr_1 \cdots dr_d.$$

We define an $L^2(\mathcal{O})$-valued cylindrical Wiener process

$$W(t) \equiv \sum_k \beta_k(t) e_k. \tag{1.13}$$

Suppose

$$\lim_{n\to\infty} m_n = \infty, \qquad \sup_n \frac{m_n^{4d}}{n} < \infty. \tag{1.14}$$

[This scaling is needed in (A.9) when verifying the exponential compact containment property (Condition 2.2) for the processes. In addition, it is also used to verify Condition 1.11(3) in Theorem 1.10.]

Let projection operator

$$P_n x \equiv \sum_{k_1=1}^{m_n} \cdots \sum_{k_d=1}^{m_n} \langle x, e_k \rangle e_k \in \operatorname{span}(e_1, \ldots, e_{m_n})$$

and for each $x \in L^2(\mathcal{O})$ fixed, we define linear operator $B(x)$ on $L^2(\mathcal{O})$ by

$$(B(x)u)(\theta) \equiv \sigma(\theta; x)u(\theta), \qquad u \in U_0 \equiv L^2(\mathcal{O}). \tag{1.15}$$

We regularize linear operator

$$B_n(x)u \equiv P_n(B(P_n x)u), \tag{1.16}$$

and arrive at an $L^2(\mathcal{O})$-valued diffusion

$$dX_n(t) = \Delta P_n X_n(t)\,dt - P_n V'(P_n X_n(t))\,dt + \frac{1}{\sqrt{n}} B_n(X_n(t))\,dW(t), \tag{1.17}$$

where the term $B_n(X_n(t))\,dW(t)$ is understood as

$$B_n(X_n(t))\,dW(t) = \sigma(\cdot; P_n X_n(t)) \sum_{k_1=1}^{m_n} \cdots \sum_{k_d=1}^{m_n} d\beta_k(t) e_k.$$

Equation (1.17) can be written in the form of (1.7). Let $\omega = \sup_{-\infty<r<\infty} |V''(r)|$ and

$$\begin{aligned} C_n x &= \Delta P_n x - P_n V'(P_n x) + P_n V'(0) - \omega P_n x, \qquad x \in L^2(\mathcal{O}), \\ F_n(x) &= -P_n V'(0) + \omega P_n x. \end{aligned} \tag{1.18}$$

Then $C_n 0 = 0$, $C_n$ is $m$-dissipative in $L^2(\mathcal{O})$ (Lemmas A.2 and A.3) and

$$dX_n(t) = C_n X_n(t)\,dt + F_n(X_n(t))\,dt + \frac{1}{\sqrt{n}} B_n(X_n(t))\,dW(t).$$

To prove the large deviation theorem, we assume:

CONDITION 1.3.

(1) $V \in C^2(R)$ and

$$\sup_{-\infty<r<\infty} |V''(r)| < \infty.$$

(2) There exist $c_1, c_2 > 0$, such that

$$V(r) \geq c_1 + c_2 r^2.$$

(3) The $\varphi(\theta, r_1, \ldots, r_k) : \mathcal{O} \times R^k \to R$ in (1.9) is bounded continuous and Lipschitz in $r_1, \ldots, r_k$, uniformly with respect to $\theta$.

Condition 1.3(3) implies that operator $B(x)$ defined by (1.15) is Lipschitz in $x$:

$$\sup_{x\neq y}\frac{|||B(x)-B(y)|||}{\|x-y\|_{L^2(\mathcal{O})}}<\infty,$$

where $|||\cdot|||$ is the operator norm

$$|||B(x)|||\equiv\sup_{\|u\|_{L^2(\mathcal{O})}\le 1}\|B(x)u\|_{L^2(\mathcal{O})}.$$

For each $n$ fixed, (1.17) can actually be represented as a finite-dimensional stochastic ODE; therefore the existence and uniqueness of the solution hold by standard finite-dimensional results.

Applying the main theorem of this paper, Theorem 1.10, we have:

THEOREM 1.4. *Under Condition 1.3 and scaling relation (1.14), the solutions $X_n(t)$ of (1.17) satisfy a large deviation principle in $C_{L^2(\mathcal{O})}[0,\infty)$ with good rate function $I$ as defined in (1.31).*

With a mild amount of additional work, assuming $\inf_{\theta,x}\sigma(\theta,x)>0$, the rate function $I$ can be represented more explicitly:

$$\begin{aligned} I(x)&=I_0(x(0))\\ &\quad+\frac{1}{2}\int_0^\infty\int_{\mathcal{O}}\left|\frac{(\partial/\partial t)x(t,\theta)-\Delta_\theta x(t,\theta)+V'(x(t,\theta))}{\sigma(\theta,x(t))}\right|^2 d\theta\,dt. \end{aligned}\tag{1.19}$$

Feng and Kurtz [18] discuss this type of representation in general. See Section 4.1 for an outline of the approach.

EXAMPLE 1.5 (Stochastic Cahn–Hilliard equation). We still consider $\mathcal{O}=[0,1)^d$ with periodic boundary condition, but with the restriction $d=1,2$ or 3 now. We consider stochastic perturbation of the Cahn–Hilliard equation formally given by

$$\frac{\partial}{\partial t}Y_n(t,\theta)=\Delta(-\Delta Y_n(t,\theta)+V'(Y_n(t,\theta)))+\frac{1}{\sqrt{n}}\frac{\partial^{d+1}}{\partial t\,\partial\theta_1\cdots\theta_d}\beta(t,\theta),\tag{1.20}$$

where $\beta(t,\theta)$ is a Brownian sheet on $[0,\infty)\times\mathcal{O}$; or, equivalently,

$$\frac{\partial}{\partial t}Y_n(t,\theta)+\mathrm{div}_\theta(\nabla(\Delta Y_n(t,\theta)-V'(Y_n(t,\theta))))=\frac{1}{\sqrt{n}}\frac{\partial^{d+1}}{\partial t\,\partial\theta_1\cdots\theta_d}\beta(t,\theta).$$

$Y_n$ is asymptotically conserved in the sense that, in the $n\to+\infty$ limit $Y$, $\int Y(t,\theta)\,d\theta$ is constant in time. As in Example 1.2, such an equation has extensive applications in material science. Motivated by the same reason as before, we only rigorously study the large deviation for a variant of (1.20), which is defined on finite Fourier modes. The number of modes goes to infinity slowly, as $n$ goes to infinity.

We assume:

CONDITION 1.6.

(1) $V \in C^3(R)$ and

$$\sup_{-\infty<r<\infty} |V''(r)| + |V'''(r)| < \infty.$$

(2) There exist $c_1, c_2 > 0$, such that

$$V(r) \geq c_1 + c_2 r^2.$$

(3)

$$\lim_{n\to\infty} m_n = \infty, \qquad \sup_n \frac{m_n^{3d}}{n} < \infty. \tag{1.21}$$

The above scaling requirement on $m_n$ is needed for reasons similar to those in the previous Allen–Cahn example. See the proof of (A.14) and the requirement in Condition 1.11(3).

We choose $E = U_0 = L^2(\mathcal{O})$. We define $e_1, \dots, e_k, \dots,$ $P_n$ and $W$ as in Example 1.2 and

$$B_n u \equiv P_n u \qquad \forall\, u \in U_0 = L^2(\mathcal{O}).$$

We consider $L^2(\mathcal{O})$-valued diffusions:

$$dX_n(t) = \Delta P_n(-\Delta P_n X_n(t)\, dt + V'(P_n X_n(t)))\, dt + \frac{1}{\sqrt{n}} B_n\, dW(t), \tag{1.22}$$

where

$$B_n dW(t) = \sum_{k_1=1}^{m_n} \cdots \sum_{k_d=1}^{m_n} e_k\, d\beta_k(t).$$

Let $\omega = \frac{1}{4}\sup_r |V''(r)|^2$,

$$C_n x \equiv \Delta P_n(-\Delta P_n x + V'(P_n x)) - \omega P_n x, \tag{1.23}$$

and

$$F_n(x) = \omega P_n x.$$

Then (1.22) can be written in the form of (1.7):

$$dX_n(t) = C_n X_n(t)\, dt + F_n(X_n(t))\, dt + \frac{1}{\sqrt{n}} B_n\, dW(t).$$

THEOREM 1.7. *Under Condition 1.6 and the scaling relation (1.21), the solutions $X_n(t)$ of (1.22) satisfy a large deviation principle in $C_{L^2(\mathcal{O})}[0,\infty)$ with good rate function $I$ as in (1.31).*

As in the stochastic Allen–Cahn example, the rate function $I$ can be further simplified:

$$I(x) = I_0(x(0)) \tag{1.24}$$
$$+ \frac{1}{2}\int_0^\infty \int_{\mathcal{O}} \left| \frac{\partial}{\partial t} x(t,\theta) - \Delta_\theta(-\Delta_\theta x(t,\theta) + V'(x(t,\theta)))\right|^2 d\theta\, dt.$$

Another type of stochastic perturbation [30] to the Cahn–Hilliard equation could also be interesting:

$$\frac{\partial}{\partial t} Y_n(t,\theta) + \mathrm{div}_\theta\left(\nabla_\theta(\Delta_\theta Y_n(t,\theta) - V'(Y_n(t,\theta))) + \frac{1}{\sqrt{n}}\frac{\beta(\partial t, \partial\theta)}{\partial t\,\partial\theta}\right) = 0,$$

where

$$\beta(t,\theta) = (\beta_1(t,\theta), \dots, \beta_d(t,\theta))$$

with each $\beta_k$ an independent real-valued space–time Brownian sheet. Large deviation for a lattice version of such an equation is considered in [18].

EXAMPLE 1.8 (Stochastic quasilinear equation with viscosity). Let $\mathcal{O} = [0,1)$ with periodic boundary. Suppose $\phi \in C^1(R)$ and $\sup_r |\phi'(r)| < \infty$. We consider the following formally defined equation:

$$\frac{\partial}{\partial t} Y_n(t,\theta) + \frac{\partial}{\partial\theta}\phi(Y_n(t,\theta)) = \alpha \frac{\partial^2}{\partial\theta^2} Y_n(t,\theta) + \frac{1}{\sqrt{n}}\frac{\partial^2}{\partial t\,\partial\theta}\beta(t,\theta),$$

where $\beta(t,\theta)$ is a Brownian sheet, $(t,\theta) \in [0,\infty) \times \mathcal{O}$.

As before, let $E = U = L^2(\mathcal{O})$ with norm $\|x\|^2 = \int_{\mathcal{O}} x^2(\theta)\, d\theta$. $\{e_1, \dots, e_k, \dots\}$ is the complete orthonormal basis as defined in (1.11) with $d = 1$. Define $L^2(\mathcal{O})$-valued cylindrical Wiener process

$$W(t) = \sum_{k=1}^\infty \beta_k(t) e_k$$

where $\{\beta_1, \beta_2, \dots\}$ are i.i.d. standard Brownian motion. Let

$$\lim_{n\to\infty} m_n = \infty, \qquad \sup_n \frac{m_n^3}{n} < \infty \tag{1.25}$$

and projection

$$P_n x = \sum_{k=1}^{2m_n} \langle x, e_k\rangle e_k.$$

We consider a regularized $L^2(\mathcal{O})$-valued diffusion equation,

$$dX_n(t) = \alpha\Delta_\theta P_n X_n(t)\,dt - P_n\,\partial_\theta\phi(P_nX_n(t))\,dt + \frac{1}{\sqrt{n}}B_n\,dW(t), \tag{1.26}$$

where $B_n = P_n$ and

$$B_n\,dW(t) = \sum_{k=1}^{2m_n} e_k d\beta_k(t).$$

The scaling on $m_n$ is needed for the same reason as in the previous two examples.

Let $\omega = (\sup_r \phi'(r))^2/(4\alpha)$,

$$C_n x = \alpha\Delta_\theta P_n x - P_n\,\partial_\theta\phi(P_n x) - \omega P_n x, \qquad x \in L^2(\mathcal{O}), \tag{1.27}$$

and

$$F_n(x) = \omega P_n x, \qquad F(x) = \omega x.$$

Then $X_n$ satisfies (1.7).

THEOREM 1.9. *Suppose that* (1.25) *holds and that* $\{X_n(0)\}$ *satisfies a large deviation principle with rate function* $I_0$. *Then* $X_n \in C_{L^2(\mathcal{O})}[0,\infty)$ *satisfy a large deviation principle with rate function given by* (1.31).

With additional work, it can be shown that the rate function $I$ admits the following form:

$$I(x) = I_0(x(0)) + \frac{1}{2}\int_0^\infty \int_{\mathcal{O}} \left|\frac{\partial}{\partial t}x - \alpha\frac{\partial^2}{\partial\theta^2}x + \frac{\partial}{\partial\theta}\phi(x)\right|^2 d\theta\,dt.$$

1.3. *Technical assumptions and main results.* Let $(E, \|\cdot\|)$ and $(U_0, \|\cdot\|_{U_0})$ be two separable real Hilbert spaces. Let $\{e_1, e_2, \dots\}$ be a complete orthonormal basis of $U_0$. We define $W$, a *cylindrical Wiener process* on $U_0$, by

$$W(t) = \sum_{k=1}^{\infty} e_k\beta_k(t), \qquad t \geq 0, \tag{1.28}$$

where $\beta_1, \beta_2, \dots$ are i.i.d. real-valued standard Brownian motions with respect to a filtration $\mathcal{F}$. $X_n(0)$ is independent of $\mathcal{F}$ and we write $\mathcal{F}^n = \mathcal{F} \vee \sigma(X_n(0))$. The infinite sum in (1.28) does not converge in $(U_0, \|\cdot\|_{U_0})$. However, we can always embed $U_0$ continuously into another separable real Hilbert space $(U_1, \|\cdot\|_{U_1})$, and as far as the embedding is Hilbert–Schmidt, the right-hand side of (1.28) converges in $(U_1, \|\cdot\|_{U_1})$. For example, let $\lambda_k > 0$ be such that $\sum_{k=1}^\infty \lambda_k^2 < \infty$; define $U_1$ to be the completion of $U_0$ under

$$\langle u, v\rangle_{U_1} = \sum_k \lambda_k^2 \langle v, e_k\rangle_{U_0}\langle u, e_k\rangle_{U_0}.$$

Then the embedding $U_0 \to U_1$ through identity map $J$ is Hilbert–Schmidt.

Throughout this paper, we will denote by $L_2(U_0, E)$ the space of Hilbert–Schmidt operators from $U_0$ into $E$ with norm

$$\begin{aligned} \|B\|^2_{L_2(U_0,E)} &\equiv \sum_k \|Be_k\|^2 \\ &= \mathrm{Tr}(B^*B) = \mathrm{Tr}(BB^*) \qquad \forall\, B \in L_2(U_0, E). \end{aligned} \tag{1.29}$$

Therefore, $J \in L_2(U_0, U_1)$. To distinguish from this space, we will denote by $L(U_0, E)$ the space of operators from $U_0$ to $E$ which are linear and bounded. For $B \in L(U_0, E)$,

$$\|B\| = \|B\|_{L(U_0,E)} = \sup_{u \in U_0, \|u\| \le 1} \|Bu\|. \tag{1.30}$$

The following stochastic integral will be used in this paper:

$$\int_0^t B(s)\, dW(s), \qquad B(s) \in L_2(U_0, E).$$

Such an integral can be defined using telescoping Riemann summation just as the usual Itô integral in finite dimensions. Although the definition of $W$ depends on $U_1$, the integral is independent of the choice of $U_1$. For details, see Chapter 4 of [9].

We identify an operator $C$ in $E$ by its graph: $C \subset E \times E$, and denote the space of continuous functions on $E$ by $C(E)$.

Our main result in the paper is the following:

THEOREM 1.10. *Let $E, U_0$ be arbitrary separable real Hilbert spaces. Suppose the following condition holds.*

CONDITION 1.11.

(1) Operator $C_n \subset E \times E$, $F_n \in C(E)$ and $B_n(x) \in L_2(U_0, E)$ are single valued and everywhere defined. That is,

$$\mathcal{D}(C_n) = \mathcal{D}(F_n) = E, \qquad \mathcal{D}(B_n(x)) = U_0 \qquad \forall\, x \in E.$$

Moreover, they are globally Lipschitz:

$$\begin{aligned} &\|C_n x - C_n y\| + \|F_n(x) - F_n(y)\| + \|B_n(x) - B_n(y)\|_{L_2(U_0,E)} \\ &\quad \le \mathrm{Constant}_n \|x - y\|, \end{aligned}$$

for every $x, y \in E$, where the $\mathrm{Constant}_n$ may depend on $n$. [See (1.29) for the definition of $\|\cdot\|_{L_2(U_0,E)}$.]

(2) $C_n$ is $m$-dissipative on $E$; $C_n 0 = 0$ for $n = 2, 3, \ldots$; and there exists a (possibly multivalued) $m$-dissipative operator $C \subset E \times E$, with $\overline{\mathcal{D}(C)} = E$ such that $C \subset \lim_{n\to\infty} C_n$, in the sense that for each $(\xi, \eta) \in C$, there exists $\xi_n \in E$ such that $\lim_n \|\xi - \xi_n\| + \|\eta - C_n \xi_n\| = 0$.

(3) Whenever $x_n \to x_0 \in E$,

$$\lim_{n\to\infty} \frac{1}{n} |||B_n(x_n)|||^2_{L_2(U_0,E)} = 0.$$

(4) For each $x \in E$, there exist $B(x) \in L(U_0, E)$ and $F \in C(E)$ satisfying

$$\sup_{x\neq y} \frac{|||B(x) - B(y)||| + \|F(x) - F(y)\|}{\|x - y\|} < \infty,$$

where $|||\cdot|||$ is the usual operator norm in (1.30). Furthermore, for each $x_n$, $p_n \in E$ and $x_n \to x_0$, $p_n \to p_0$, we have

$$F_n(x_n) \to F(x_0) \quad \text{and} \quad \|B_n^*(x_n)p_n\|_{U_0} \to \|B^*(x_0)p_0\|_{U_0}.$$

We also assume that $X_n$ is the solution to (1.7), and that $\{X_n : n = 1, 2, \ldots\}$ satisfies the following exponential compact containment condition:

CONDITION 1.12. For each compact $K \subset E$, $T > 0$ and $a > 0$, there exists another compact set $K_{a,T} \subset E$ such that

$$\limsup_{n\to\infty} \sup_{x\in K} \frac{1}{n} \log P(X_n(t) \notin K_{a,T}, \exists\, 0 < t \leq T | X_n(0) = x) \leq -a.$$

Finally, we assume that $\{X_n(0) : n = 1, 2, \ldots\}$ satisfies the large deviation principle with good rate function $I_0$ on $E$.

Then:

(a) $\{X_n\}$ is exponentially tight;

(b) the following limit exists and defines an operator semigroup on $C_b(E)$ (the space of bounded continuous functions on $E$):

$$V(t)f(x) = \lim_{n\to\infty, y\to x} \frac{1}{n} \log E[e^{nf(X_n(t))} | X_n(0) = y];$$

(c) the large deviation principle holds for $\{X_n\}$ with good rate function $I$:

$$I(x) = I_0(x(0)) + \sup_{0\leq t_1\leq\cdots\leq t_m} \left( \sum_{i=1}^m I_{t_i - t_{i-1}}(x(t_i)|x(t_{i-1})) \right), \tag{1.31}$$

where

$$I_t(y|x) = \sup_{f\in C_b(E)} (f(y) - V(t)f(x)).$$

Theorems 1.4, 1.7 and 1.9 are all special cases of this theorem.

1.4. *Relation to other large deviation results in literature.* The term $C_n$ in (1.7) and its limit $C$ in Theorem 1.10 are allowed to be totally nonlinear. This is different than what is available in literature [4, 9, 24, 27, 29], where $C$ is restricted to be semilinear. However, this paper does not pursue generalities in the term $B_n$, as some of the above mentioned papers do.

If $C_n$ is semilinear, a good deal is known about the solution for (1.7). See [9]. However, if $C_n$ is just $m$-dissipative, very little is known for the equation. By assuming $C_n$ is Lipschitz and everywhere defined for each fixed $n$, we greatly simplified the situation. Such assumption is motivated by Examples 1.2 and 1.8, and by the fact that Yosida approximation of $m$-dissipative operators satisfies the above requirements.

In this paper (1.7) is driven by a Brownian noise $W$, which is responsible for the quadratic nonlinear term in $H_0, H_1$ (or the $\hat{H}_0, \hat{H}_1$). But in the proof of comparison principle, we actually allow much more general nonlinearity (Theorem 5.2). Therefore, it is possible that the method here can be applied to cases where the $W$ is replaced by spatial Poisson noise. We expect exponential nonlinearity in the $H_i, \hat{H}_i$'s in these cases.

1.5. *Notation.* We will frequently use the following class of test functions for localization purpose:

$$\begin{aligned}\mathcal{T} = \{\varphi \in C^2([0,\infty)) : \varphi \geq 0,\ &\text{is nondecreasing,} \\ &\text{and } \varphi(r) = \varphi(+\infty) \text{ for } r \text{ large enough}\}.\end{aligned} \tag{1.32}$$

Throughout the paper, $(E, r)$ and $(U, r_U)$ are complete separable metric spaces. Let $f$ be a function on $E$. $C(E)$ denotes continuous functions on $E$; $C_b(E)$, bounded continuous functions; $B(E)$, bounded Borel measurable functions; $M(E)$, Borel measurable functions; and $\mathcal{P}(E)$, probability measures on $E$. For $f \in M(E)$, we define $f^*$ and $f_*$ to be, respectively, the upper and lower semicontinuous smoothing of $f$:

$$f^*(x) = \limsup_{y \to x} f(y), \qquad f_*(x) = \liminf_{y \to x} f(y).$$

Let $\mathcal{O} \subset R^d$:

$$H^1(\mathcal{O}) \equiv \left\{ x(\theta) \in L^2(\mathcal{O}) : \int_{\mathcal{O}} |x(\theta)|^2 + |\nabla x(\theta)|^2 \, d\theta < \infty \right\}$$

and

$$H^2(\mathcal{O}) \equiv \left\{ x(\theta) \in L^2(\mathcal{O}) : \sum_{i,j,k=1}^{d} \int_{\mathcal{O}} |x(\theta)|^2 + |\partial_k x(\theta)|^2 + |\partial_{i,j}^2 x(\theta)|^2 \, d\theta < \infty \right\}.$$

Throughout, $(E, \|\cdot\|)$ is a real separable Hilbert space with its dual identified as itself $E^* = E$. $C^k(E)$ denotes the set of $k$th-order Fréchet differentiable functions on $E$ with continuous $k$th-order derivative; we identify the $k$th-order derivative as a $k$th-order multilinear symmetric functional. For example, by $Df(x)$, we mean the gradient of $f$ evaluated at $x$, which is identified as an element of $E^* \equiv E$ through the Taylor expansion:

$$f(x+y) = f(x) + \langle Df(x), y\rangle + \frac{1}{2}D^2 f(x)yy + o(\|y\|^2), \qquad y \in E. \tag{1.33}$$

$D^2 f(x)yz$ means a functional which is bilinear in both the $y$ and the $z$ arguments. Let $x, y \in E$; by $x \otimes y$, we mean a bounded linear operator on $E$

$$(x \otimes y)z \equiv x\langle y, z\rangle.$$

$|\cdot|$ will be used to denote either an absolute value of a number $|a|$ or the Euclidean norm of a vector in $R^d$: $|(\theta_1, \ldots, \theta_d)|^2 = \sum_{k=1}^d \theta_k^2$.

We denote the range of a generic operator $A$ in a Banach space by $\mathcal{R}(A)$ and its domain by $\mathcal{D}(A)$. We often identify an operator with its graph. $\bar{A}$ denotes the closure of the graph under the norm of the Banach space. Let $E$ be a Hilbert space. A possibly multivalued nonlinear operator $C \subset E \times E$ is said to be *m-dissipative* if and only if it is *dissipative*:

$$\langle x_1 - x_2, y_1 - y_2\rangle \le 0 \qquad \forall\, (x_i, y_i) \in C$$

and

$$\mathcal{R}(I - \alpha C) = E \qquad \forall\, \alpha > 0.$$

If, in addition, $\overline{\mathcal{D}(C)} = E$, then $C$ generates a strongly continuous contraction semigroup $S(t)$ on $E$. The following test function on $E$ is introduced to record trajectory properties of $S(t)y$. It plays a major role in the analysis of certain Hamilton–Jacobi equations (Section 5).

Definition 1.13. Let $C$ be an $m$-dissipative operator on $E$ generating a strongly continuous semigroup $S(t): t \ge 0$. The *Tataru distance function* $d_C$ is

$$d_C(x, y) \equiv \inf\{t + \|x - S(t)y\| : t \ge 0\} \qquad \forall\, x, y \in E.$$

$d_C(x, y)$ is Lipschitz ((28) on page 62 of [8]):

$$|d_C(x, y) - d_C(\hat{x}, \hat{y})| \le \|x - \hat{x}\| + \|y - \hat{y}\|.$$

Let $E$ be a general metric space again and let $H_0, H_1 \subset C_b(E) \times B(E)$ be (possibly multivalued) operators, $h \in C_b(E)$ and $\alpha > 0$. We define viscosity solutions for

$$(I - \alpha H_0)f = h \tag{1.34}$$

and

$$(I - \alpha H_1)f = h. \tag{1.35}$$

DEFINITION 1.14 (Viscosity solution).

(a) $\overline{f}$ is a *viscosity subsolution* of (1.34) if and only if $\overline{f}$ is bounded, upper semicontinuous, and for each $(f_0, g_0) \in H_0$, there exists an $\{x_n\} \subset E$ satisfying

$$\lim_{n\to\infty} (\overline{f} - f_0)(x_n) = \sup_{x\in E} (\overline{f} - f_0)(x) \tag{1.36}$$

and

$$\limsup_{n\to\infty} (\alpha^{-1}(\overline{f} - h) - (g_0)^*)(x_n) \le 0. \tag{1.37}$$

(b) $\underline{f}$ is a *viscosity supersolution* of (1.35) if and only if $\underline{f}$ is bounded, lower semicontinuous, and for each $(f_0, g_0) \in H_1$, there exists an $\{x_n\} \subset E$ satisfying

$$\lim_{n\to\infty} (f_0 - \underline{f})(x_n) = \sup_{x\in E} (f - f_0)(x) \tag{1.38}$$

and

$$\liminf_{n\to\infty} (\alpha^{-1}(\underline{f} - h) - (g_0)_*)(y_n) \ge 0. \tag{1.39}$$

DEFINITION 1.15. We say a comparison principle holds for viscosity subsolution of (1.34) and supersolution of (1.35) if

$$\overline{f} \le \underline{f},$$

for every subsolution $\overline{f}$ of (1.34) and supersolution $\underline{f}$ of (1.35).

Allowing $H_0, H_1 \subset C(E) \times M(E)$, Tataru [31, 32] and Crandall and Lions [8] define viscosity solution in a different manner. Definition 1.16 is an adaptation of their definitions when the domain of operator is chosen properly. We will explore the connection between such definition and the more general Definition 1.14, for equations arising in our large deviation context. See Section 3.4.

DEFINITION 1.16 (Tataru–Crandall–Lions).

(a) We say that $\overline{f}$ is a subsolution of (1.34), if $\overline{f}$ is bounded upper semicontinuous on $E$, and for each $x_0 \in E$ and $f_0 \in \mathcal{D}(H_0)$ satisfying

$$(\overline{f} - f_0)(x_0) = \sup_{x\in E} (\overline{f} - f_0)(x), \tag{1.40}$$

we have

$$\alpha^{-1}(\overline{f} - h)(x_0) \le (H_0 f_0)^*(x_0).$$

(b) We say that $\underline{f}$ is a supersolution of (1.35), if $\underline{f}$ is bounded lower semi-continuous on $E$, and for each $x_0 \in E$ and $f_0 \in \mathcal{D}(H_1)$ satisfying

$$(f_0 - \underline{f})(x_0) = \sup_{x \in E}(f_0 - \underline{f})(x), \tag{1.41}$$

we have

$$\alpha^{-1}(\underline{f} - h)(x_0) \geq (H_1 f_0)_*(x_0).$$

DEFINITION 1.17 (Exponential tightness and large deviation principle). Let $S$ be a complete separable metric space and let $\{X_n\}$ be $S$-valued random variables. $\{X_n\}$ is said to be *exponentially tight* if for every $a > 0$, there exists compact $\mathcal{K}_a \subset S$ such that

$$\limsup_{n \to \infty} \frac{1}{n} \log P(X_n \notin \mathcal{K}_a) < -a.$$

$\{X_n\}$ is said to satisfy the *large deviation principle* if there exists a lower semicontinuous function $I : S \to [0, +\infty]$ such that for every open set $A \subset S$,

$$-\inf_{x \in A} I(x) \leq \liminf_n \frac{1}{n} \log P(X_n \in A)$$

and for every closed set $B \subset S$,

$$\limsup_n \frac{1}{n} \log P(X_n \in B) \leq -\inf_{x \in B} I(x).$$

$I$ is called the *rate function* and it is *good* if each level set is compact.

Let $E$ be a complete separable metric space. For each $n$, let stochastic process $X_n$ have state space $E$ and let its trajectory be continuous in time. By large deviation (resp. exponential tightness) for the processes $\{X_n\}$, we apply the above definition with $S = C_E[0, \infty)$.

DEFINITION 1.18. Let $(E, q)$ be a metric space. $D \subset C_b(E)$ is said to *approximate the metric* $q$ if for each compact $K \subset E$ and $z \in K$, there exists $f_n \in D$ such that $\lim_{n \to \infty} \sup_{x \in K} |f_n(x) - q(x, z)| = 0$.

**2. A general large deviation theorem.** This section presents a general theorem which is the basis for the large deviation method in this paper. The heuristics have been explained in Section 1.1.

Let $(E, r)$ be a complete separable metric space and let $\{X_n\}$ be a sequence of $E$-valued processes with trajectories in $C_E[0, \infty)$. Suppose $A_n \subset B(E) \times B(E)$ is possibly multivalued. Let $X_n$ be a solution to the $A_n$-*martingale problem*. That is, there is a filtration $\mathcal{F}_t^n$, such that

$$f(X_n(t)) - f(X_n(0)) - \int_0^t g(X_n(s))\, ds \qquad \forall\, (f, g) \in A_n$$

is a martingale. We will work under the following regularity condition.

Condition 2.1. For each $n = 2, 3, \ldots$, let $A_n \subset B(E) \times B(E)$. We assume existence and uniqueness hold for the martingale problem for $A_n$ with $X_n \in C_E[0, +\infty)$ for each initial distribution $\mu \in \mathcal{P}(E)$. Let $P_x^n \in \mathcal{P}(C_E[0, \infty))$ denote the distribution of the solution of the martingale problem for $A_n$ with $X_n(0) = x \in E$; we assume that the mapping $x \to P_x^n$ is Borel measurable taking the weak topology on $\mathcal{P}(C_E[0, \infty))$ (cf. Theorem 4.4.6 of [15]).

Define $H_n \subset B(E_n) \times B(E_n)$ by

$$H_n f = \frac{1}{n} e^{-nf} A_n e^{nf}, \qquad e^{nf} \in \mathcal{D}(A_n),$$

or if $A_n$ is multivalued,

$$H_n = \left\{ \left( f, \frac{1}{n} e^{-nf} g \right) : (e^{nf}, g) \in A_n \right\}.$$

The following is an exponential version of the *uniform compact containment* condition in the weak convergence theory.

Condition 2.2. For each compact $K \subset E$, $T > 0$ and $a > 0$, there exists a compact $K_{a,T} \subset E$ such that

$$\limsup_{n\to\infty} \sup_{x \in K} \frac{1}{n} \log P(X_n(t) \notin K_{a,T}, \text{ for some } 0 \le t \le T | X_n(0) = x) \le -a.$$

The following is an adaptation of Theorem 7.18 of [18]. In the adaptation, we also used a result of exponential tightness (Corollary 4.19), a variant of the Stone–Weierstrass theorem (Lemma A.8) and a technical estimate (Lemma 7.19). All the reference labels refer to [18].

Theorem 2.3. *Let Condition 2.1 be satisfied. In addition, we assume the following:*

(1) *Convergence of generators. There exist $H_0, H_1 \subset C_b(E) \times B(E)$ which are limits of the $H_n$'s in the following sense:*
   (a) *For each $(f, g) \in H_0$, there exist some $(f_n, g_n) \in H_n$ such that*

$$\sup_n \left( \sup_x |f_n(x)| + \sup_x |g_n(x)| \right) < \infty,$$

*and that for each $x_n \to x_0$, we have*

$$\lim_{n\to\infty} f_n(x_n) = f(x_0), \qquad \limsup_{n\to\infty} g_n(x_n) \le g^*(x_0).$$

(b) *For each* $(f,g) \in H_1$, *there exist some* $(f_n, g_n) \in H_n$ [*possibly different than those in* (a)] *such that*

$$\sup_n \left( \sup_x |f_n(x)| + \sup_x |g_n(x)| \right) < \infty,$$

*and that for each* $x_n \to x_0$, *we have*

$$\lim_{n\to\infty} f_n(x_n) = f(x_0), \qquad g_*(x_0) \le \liminf_{n\to\infty} g_n(x_n).$$

*There exist* $F \subset C_b(E)$ *which approximate the metric* $q \equiv r \wedge 1$ *(Definition 1.18), and for each* $f \in F$ *and* $\lambda > 0$, $\lambda f \in \mathcal{D}(H_0)$.

(2) *Uniform exponential compact containment. Condition 2.2 holds.*

(3) *Comparison principle. There exist a subset* $D \subset C_b(E)$ *and* $\alpha_0 > 0$, *such that for each* $h \in D$ *and* $0 < \alpha < \alpha_0$, *the comparison principle (Definition 1.15) holds for subsolution (in the sense of Definition 1.14) of*

$$(I - \alpha H_0) f = h,$$

*and supersolution (Definition 1.14) of*

$$(I - \alpha H_1) f = h.$$

$D$ *contains an algebra that separates points and vanishes nowhere* [*i.e., for each* $x \in E$, *there exists* $f$ *belonging to this algebra such that* $f(x) \neq 0$].

*Define* $\{V_n(t)\}$ *on* $B(E)$ *by*

$$V_n(t) f(x) = \frac{1}{n} \log E[e^{n f(X_n(t))} | X_n(0) = x].$$

*If* $\{X_n(0)\}$ *satisfies a large deviation principle in* $E$ *with a good rate function* $I_0$, *then:*

(a) *limit*

$$V(t) f(x) \equiv \lim_{n\to\infty} V_n(t_n) f(x_n)$$

*exists for every* $f \in C_b(E)$, $t_n \to t$, *and* $x_n \to x$. $V(t)$ *forms a nonlinear semigroup on* $C_b(E)$.

(b) *for each* $0 \le t_1 \le \cdots \le t_k < \infty$, $\{(X_n(t_1), \ldots, X_n(t_k)) : n = 1, 2, \ldots\}$ *is exponentially tight in* $E^k$ *and satisfies the large deviation principle with good rate function*

$$\begin{aligned} I_{t_1,\ldots,t_k}&(x_1,\ldots,x_k) \\ &= \sup_{f_1,\ldots,f_k \in D} \{f_1(x_1) + \cdots + f(x_k) \\ &\qquad - \Lambda_0(V(t_1)(f_1 + V(t_2 - t_1) \end{aligned} \tag{2.1}$$

$$\times (f_2 + \cdots + V(t_k - t_{k-1}) f_k) \ldots))\}$$

$$= \inf_{x_0 \in E} \left\{ I_0(x_0) + \sum_{i=1}^{k} I_{t_i - t_{i-1}}(x_i | x_{i-1}) \right\},$$

*where*

$$\Lambda_0(f) = \lim_{n \to +\infty} \frac{1}{n} \log E[e^{nf(X_n(0))}] \qquad \forall f \in C_b(E)$$

[*the limit exists by* (1.1)], *and*

$$I_t(y|x) = \sup_{f \in C_b(E)} (f(y) - V(t)f(x)). \tag{2.2}$$

(c) $\{X_n\}$ *is exponentially tight in* $C_E[0, \infty)$ *and satisfies the large deviation principle with good rate function:*

$$\begin{aligned} I(x) &= \sup_{k=1,2,\ldots} I_{t_1,\ldots,t_k}(x(t_1), \ldots, x(t_k)) \\ &= \sup_{k=1,2,\ldots} \sup_{0<t_1<t_2<\cdots<t_k} \left( I_0(x(0)) + \sum_{i=1}^{k} I_{t_i - t_{i-1}}(x(t_i)|x(t_{i-1})) \right). \end{aligned} \tag{2.3}$$

REMARK 2.4. In view of the duality in (1.2), the form of rate function in the first identity of (2.1) should be expected. The second equality follows by the Markovian property of the $X_n$'s. The rate function (2.3) follows from the finite-dimensional large deviation result in (b), and from a well-known projective limit argument [10, 11, 28].

Large deviation behavior of the $\{X_n\}$ and the exponential tightness imply that the rate function $I$ is good. That is, $I$ has compact level sets. See part (b) of Lemma 1.2.18 of [12]. Similarly, $I_{t_1,\ldots,t_k}$ has compact level sets in $E^k$.

Theorem 2.3 can be applied to situations other than small perturbation type problems; we refer the reader to [18] for further examples.

In the rest of this paper we apply the above theorem to the general problem considered in Theorem 1.10. Step 1 is verified in Section 3.3; see Lemma 3.4. The condition in step 2 is assumed in Theorem 1.10 as Condition 1.12. It is verified for Examples 1.2, 1.5 and 1.8 in Section A.2. The comparison principle in step 3 is stated in Lemma 3.10, with details of the actual proof carried out in Section 5.

Before closing this section, we illustrate how the classical Freidlin–Wentzell theory follows from Theorem 2.3.

EXAMPLE 2.5 (The Freidlin–Wentzell theory). Let $E = R^d$, and let $W$ be a standard $d$-dimensional Brownian motion. Assume that $b_j, \sigma_{ij} \in C_b(R^d)$ are Lipschitz continuous for $i, j = 1, \ldots, d$. We denote $b(x) = (b_1(x), \ldots, b_d(x)) : R^d \to R^d$, and let $d \times d$-matrix $\sigma(x) = (\sigma_{ij}(x))$. Let $X_n$ be the solution to

$$dX_n(t) = b(X_n(t))\,dt + \frac{1}{\sqrt{n}}\sigma(X_n(t))\,dW(t).$$

This is a special case of (1.7).

Let $D = \{f : f = f_0 + c, f_0 \in C_0^2(R^d), c \in R\}$ where $C_0^2(R^d)$ is the collection of functions with compact support and with continuous derivative up to the second order. We denote by $D^2 f(x) = (\partial_{ij}^2 f(x))_{i,j}$ the Hessian matrix of $f$. By Itô's formula, if we take

$$A_n f(x) = b(x)\nabla f(x) + \frac{1}{2n}\operatorname{Tr}(D^2 f(x)\sigma(x)\sigma^T(x)), \qquad f \in D,$$

then Condition 2.1 is satisfied. The transformed generator

$$H_n f(x) = b(x)\nabla f(x) + \frac{1}{2}|\sigma^T(x)\nabla f(x)|^2 + \frac{1}{2n}\operatorname{Tr}(D^2 f(x)\sigma(x)\sigma^T(x)), \qquad f \in D.$$

If we let $H_0 = H_1 = H$ with

$$Hf(x) = b(x)\nabla f(x) + \tfrac{1}{2}|\sigma^T(x)\nabla f(x)|^2, \qquad f \in D,$$

then the convergence conditions in part (1) of Theorem 2.3 are satisfied.

The assumptions on $\sigma, b$ imply that they grow at most linearly:

$$\sum_{i,j,k=1}^{d}(|\sigma_{ij}(x)| + |b_k(x)|) \le c_1 + c_2|x|.$$

One can use such estimate to verify the uniform exponential compact containment condition in Theorem 2.3. This is shown in Example 4.23 of [18] using a stochastic Lyapunov function technique.

Let $h \in C_b(R^d)$ and $\alpha > 0$; the comparison principle for

$$(I - \alpha H)f = h \tag{2.4}$$

follows from results in [5]. Details on its proof can also be found in Chapters 9.4 and 10.3 of [18].

Consequently, by Theorem 2.3, the large deviation principle holds for $\{X_n\}$.

Let $\int_0^T |u(s)|^2\,ds < +\infty$ for each $T > 0$ and consider

$$\dot{x}(t) = b(x(t)) + \sigma(x(t))u(t). \tag{2.5}$$

We define

$$R_\alpha h(x_0) = \sup\left\{\int_0^\infty e^{-\alpha^{-1}s}(\alpha^{-1}h(x(s)) - \tfrac{1}{2}|u(s)|^2)\,ds : (x, u) \text{ satisfies (2.5)}, \; x(0) = x_0\right\}.$$

Then it can be shown that $R_\alpha h \in C_b(R^d)$, it is the unique solution to (2.4), and

$$V(t)h(x_0) = \lim_{k\to+\infty} R^k_{t/k} h(x_0)$$
$$= \sup\left\{\int_0^t \tfrac{1}{2}|u(s)|^2\,ds + h(x(t)) : \dot{x} = b(x) + \sigma(x)u, x(0) = x_0\right\}.$$

All these can be rigorously justified using the dynamic programming principle. Suppose that $\sigma^{-1}(x)$ exists. Plug the above expression on $V$ in (2.3) and (2.2); we obtain the simplified representation

$$I(x) = \int_0^\infty |\sigma^{-1}(x)(\dot{x} - b(x))|^2\,ds.$$

Such result is known as the Freidlin–Wentzell theory [23].

All the above claims are well-known results in control theory and first-order Hamilton–Jacobi equation literature. In [18], rigorous proofs are provided and summarized again.

**3. Large deviation for diffusions in Hilbert space.** Throughout this section we assume Condition 1.11 holds and $X_n$ is a solution of (1.7). Both $E$ and $U_0$ are real separable Hilbert spaces.

3.1. *Semigroup on Hilbert spaces.* We first recall some basic facts about semigroup on $E$. These facts will be used later in the paper.

We assumed, in Condition 1.11(2), that $C$ is $m$-dissipative on $E$, and $\overline{\mathcal{D}(C)} = E$. By Crandall and Liggett's [6] semigroup generation theorem, it generates a strongly continuous contraction semigroup on $E$:

$$S(t)x = \lim_{n\to\infty}\left(I - \frac{t}{n}C\right)^{-n} x \qquad \forall x \in E,$$

and

$$\|S(t)x - S(t)y\| \le \|x - y\| \qquad \forall t > 0, x, y \in E.$$

Since $(0,0) \in C$ (or $C0 = 0$ if $C$ is single valued), $0 = (I - \alpha C)^{-1}0$, $S(t)0 = 0$ and $\|S(t)x\| \le \|x\|$.

Definition 3.1 (Canonical restriction of $C$). We denote

$$\|Cx\| \equiv \inf\{\|y\| : (x,y) \in C\} \qquad \forall x \in \mathcal{D}(C)$$

and define a single-valued $C^0 \subset E \times E$, called *the canonical restriction of* $C$, by

$$C^0 x = \{z : (x,z) \in C, \|z\| = \|Cx\|\}.$$

Then the following holds.

LEMMA 3.2.

(1) *$\mathcal{D}(C^0) = \mathcal{D}(C)$ and $C^0$ is single valued (Lemma 2.19 of [26]).*
(2) *$C^0$ is the infinitesimal generator of $S(t)$ in the sense that*

$$C^0 x = \lim_{h \to 0+} \frac{1}{h}(S(h)x - x), \qquad x \in \mathcal{D}(C) \tag{3.1}$$

*(Corollary 4.19 in [26]).*
(3) *Let $f \in C^1(E)$; then*

$$\langle Df(\xi), C^0 \xi \rangle = \lim_{r \to 0+} \frac{f(S(r)\xi) - f(\xi)}{r} \qquad \forall \xi \in \mathcal{D}(C). \tag{3.2}$$

DEFINITION 3.3 [Directional derivative along the trajectory of $S(t)$]. Suppose $f \in C(E)$ is Lipschitz continuous. We define

$$D_C^+ f(x) = \limsup_{h \to 0+, y \to x} \frac{1}{h}(f(S(h)y) - f(y))$$

and

$$D_C^- f(x) = \liminf_{h \to 0+, y \to x} \frac{1}{h}(f(S(h)y) - f(y)).$$

$D_C^+ f : E \to [-\infty, +\infty]$ is upper semicontinuous, and $D_C^- f : E \to [-\infty, +\infty]$ is lower semicontinuous (Lemma 2.3 in [8]).

We list two useful properties of the Tataru distance function $d_C$ (Definition 1.13):

$$d_C(x, y) - d_C(\hat{x}, \hat{y}) \le \|x - \hat{x}\| + \|y - \hat{y}\| \tag{3.3}$$

and

$$\frac{d_C(S(r)x, y) - d_C(x, y)}{r} \le 1, \qquad D_C^+ d_C(\cdot, y) \le 1. \tag{3.4}$$

See page 62 of [8] for proof.

3.2. *The martingale problem.* Recall that $E$ is a real separable Hilbert space, and that $C_n$ is single valued, everywhere defined and Lipschitz on the $E$ ($C_n$ is usually some regularization of the $C$). Let $f \in C^2(E)$ be such that $Df(x) = 0$ when $\|x\|$ is sufficiently large; we define linear operator $A_n \subset C_b(E) \times B(E)$ by

$$A_n f(x) = \langle Df(x), C_n x + F_n(x) \rangle + \frac{1}{2n} \mathrm{Tr}[D^2 f(x) B_n(x) B_n^*(x)]. \tag{3.5}$$

By Condition 1.11, an infinite-dimensional version of the Itô formula applies (e.g., Theorem 4.17 of [9]). In addition, Condition 1.11(1) is a strong enough assumption so that the results in Chapter 9 of [9] (regarding Markov property and regularity for the initial conditions) apply. Therefore Condition 2.1 is satisfied.

We next compute nonlinear operator $H_n \subset C_b(E) \times B(E)$ by

$$\begin{aligned}
H_n f(x) &\equiv \frac{1}{n} e^{-nf} A_n e^{nf}(x) \\
&= \frac{1}{n} e^{-nf}(x) \langle D e^{nf}(x), C_n x + F_n(x) \rangle \\
&\quad + \frac{1}{2n^2} e^{-nf(x)} \operatorname{Tr}[D^2 e^{nf}(x) B_n(x) B_n^*(x)] \\
&= \langle Df(x), C_n x + F_n(x) \rangle \qquad (3.6) \\
&\quad + \frac{1}{2n^2} \operatorname{Tr}[(D(nf)(x) \otimes D(nf)(x) + D^2(nf)(x)) B_n(x) B_n^*(x)] \\
&= \langle Df(x), C_n x + F_n(x) \rangle \\
&\quad + \frac{1}{2} \| B_n^*(x) Df(x) \|_{U_0}^2 + \frac{1}{2n} \operatorname{Tr}[D^2 f(x) B_n(x) B_n^*(x)],
\end{aligned}$$

where $(x \otimes y)z \equiv x\langle y, z \rangle$. The last step above needs some justification: let $\{\hat{e}_1, \ldots, \hat{e}_k, \ldots\}$ be a complete orthonormal basis of $E$. Then

$$\begin{aligned}
&\frac{1}{2n^2} \operatorname{Tr}[(D(nf)(x) \otimes D(nf)(x)) B_n(x) B_n^*(x)] \\
&\quad = \frac{1}{2} \sum_{k=1}^{\infty} \langle (Df(x) \otimes Df(x)) B_n(x) B_n^*(x) \hat{e}_k, \hat{e}_k \rangle \\
&\quad = \frac{1}{2} \sum_{k=1}^{\infty} \langle Df(x), \hat{e}_k \rangle \langle Df(x), B_n(x) B_n^*(x) \hat{e}_k \rangle \\
&\quad = \frac{1}{2} \langle Df(x), B_n(x) B_n^*(x) Df(x) \rangle \\
&\quad = \frac{1}{2} \| B_n^*(x) Df(x) \|_{U_0}^2.
\end{aligned}$$

3.3. *Convergence of the $H_n$'s.* Formally, we expect the limit of $H_n f$ to be given by

$$Hf(x) = \langle Df(x), Cx + F(x) \rangle + \tfrac{1}{2} \| B^*(x) Df(x) \|_{U_0}^2.$$

However, the above does not make sense for $x \notin \mathcal{D}(C)$. As commented in Section 1.1, we have to replace $H$ by $H_0, H_1$; then by selecting test functions $f$ carefully, we can estimate the limit from above by $H_0 f$ and from below by $H_1 f$. The class of test functions has to be large enough so that the comparison principle (Sections 3.4 and 5) can be proved.

This is what we will carry out rigorously next.

By Condition 1.11, $C_n$ and $C$ generate, respectively, strongly continuous contraction semigroup $S_n(t)$ and $S(t)$ on $E$, $S_n(t)0 = 0$, $S(t)0 = 0$ and $\|S_n(t)x\| \le \|x\|$, $\|S(t)x\| \le \|x\|$.

We derive the limit operators $H_0, H_1$ in Theorem 2.3 through several steps.

First, recall definitions of the *canonical restriction of* $C$ in Definition 3.1 and of the *Tataru distance function* $d_C$ in Definition 1.13. $d_C$ is Lipschitz; however, it may not be differentiable in $x$. We introduce smooth approximations of it first.

By Condition 1.11, both $C_n$ and $C$ are $m$-dissipative, and $C \subset \lim_n C_n$. By the Crandall–Liggett semigroup convergence theorem ([6]; see also Theorem 6.8 of [26]),

$$\lim_{n\to\infty} \sup_{0\le t\le T} \|S_n(t)y - S(t)y\| = 0 \qquad \forall\, y \in E, T > 0.$$

Let $\lim_{n\to\infty} a_n = \infty$; we define

$$\phi_\varepsilon(r) = \left(\sqrt{\varepsilon} + \frac{r-\varepsilon}{2\sqrt{\varepsilon}} - \frac{(r-\varepsilon)^2}{8\varepsilon\sqrt{\varepsilon}}\right) I(0 \le r < \varepsilon) + \sqrt{r} I(r \ge \varepsilon), \tag{3.7}$$

$$h_{\varepsilon,y}(x) \equiv \inf_{t\ge 0}\{t + \phi_\varepsilon(\|x - S(t)y\|^2)\}, \tag{3.8}$$

$$h_{n,\varepsilon,y}(x) \equiv -\frac{1}{a_n} \log \int_0^\infty e^{-a_n\{t+\phi_\varepsilon(\|x-S_n(t)y\|^2)\}}\, dt. \tag{3.9}$$

Then by Lemma A.12,

$$\lim_{\varepsilon\to 0+} \sup_{x\in E} |h_{\varepsilon,y}(x) - d_C(x,y)| = 0$$

and

$$\lim_{n\to\infty} \sup_{x\in K} |h_{\varepsilon,y}(x) - h_{n,\varepsilon,y}(x)| = 0$$

for each compact $K \subset E$. Later, we may drop the $y$ in the subindex if no confusion can occur.

Recall the definition of $\mathcal{T}$ in (1.32). We now define $H_0$ and $H_1$:

(a) Let

$$\mathcal{D}(H_0) = \{f(x) : f(x) = \varphi_1(\|x-\xi\|^2) + \varphi_2(h_{\varepsilon,y_1}(x)) + \cdots + \varphi_{k+1}(h_{\varepsilon,y_k}(x)),$$
$$\forall\, \varphi_i \in \mathcal{T}, \xi \in \mathcal{D}(C), y_j \in E, k = 1, 2, \dots\}.$$

For $g(x)=\varphi_1(\|x-\xi\|^2)$ and $f(x)=g(x)+\varphi_2(h_{\varepsilon,y_1}(x))+\cdots+\varphi_{k+1}(h_{\varepsilon,y_k}(x))\in\mathcal{D}(H_0)$, we define

$$\begin{aligned}H_0f(x)&=2\varphi_1'(\|x-\xi\|^2)\langle x-\xi,C^0\xi\rangle\\&\qquad+\Bigl(\sup_{r\ge0}\varphi_2'(r)+\cdots+\sup_{r\ge0}\varphi_{k+1}'(r)\Bigr)\\&\qquad+\sup_{\|q\|\le\varphi_2'(h_{\varepsilon,y_1}(x))+\cdots+\varphi_{k+1}'(h_{\varepsilon,y_k}(x))}(\langle F(x),Dg(x)+q\rangle\\&\qquad+\tfrac12\|B^*(x)(Dg(x)+q)\|_{U_0}^2).\end{aligned}\tag{3.10}$$

By item (1) of Lemma A.12, and the fact that $Dg(x)=0$ when $\|x\|$ is sufficiently large, we have

$$\sup_{x\in E}\sup_{\|q\|\le\varphi_2'(h_{\varepsilon,y_1}(x))+\cdots+\varphi_{k+1}'(h_{\varepsilon,y_k}(x))}\|\langle F(x),Dg(x)+q\rangle\|+\tfrac12\|B^*(x)(Dg(x)+q)\|_{U_0}^2<\infty.$$

Consequently, $H_0\subset C_b(E)\times B(E)$.

(b) Let

$$\begin{aligned}\mathcal{D}(H_1)=\{f(x):f(x)&=-\varphi_1(\|x-\xi\|^2)\\&\quad-\varphi_2(h_{\varepsilon,y_1}(x))-\cdots-\varphi_{k+1}(h_{\varepsilon,y_k}(x)),\\&\quad\forall\varphi_i\in\mathcal{T},\xi\in\mathcal{D}(C),y_j\in E,k=1,2,\dots\}.\end{aligned}$$

Let $f(x)=g(x)-\varphi_2(h_{\varepsilon,y_1}(x))-\cdots-\varphi_{k+1}(h_{\varepsilon,y_k}(x))\in\mathcal{D}(H_1)$, where $g(x)=-\varphi_1(\|x-\xi\|^2)$. We define

$$\begin{aligned}H_1f(x)&=-2\varphi_1'(\|x-\xi\|^2)\langle x-\xi,C^0\xi\rangle\\&\qquad-\Bigl(\sup_{r\ge0}\varphi_2'(r)+\cdots+\sup_{r\ge0}\varphi_{k+1}'(r)\Bigr)\\&\qquad+\inf_{\|q\|\le\varphi_2'(h_{\varepsilon,y_1}(x))+\cdots+\varphi_{k+1}'(h_{\varepsilon,y_2}(x))}(\langle F(x),Dg(x)+q\rangle\\&\qquad+\tfrac12\|B^*(x)(Dg(x)+q)\|_{U_0}^2).\end{aligned}\tag{3.11}$$

The reason for using $\varphi_i$ is to localize the test function $f$ and the $H_0f,H_1f$ so that $H_0,H_1\subset C_b(E)\times B(E)$, a condition required by Theorem 2.3. This is the main reason that the two operators have such complicated forms, instead of the simpler forms used by Crandall and Lions [8] using the $D_C^+$ in Definition 3.3. We note that $H_0$ and $H_1$ are both single valued.

Let

$$F=\{f(x)=\varphi_1(\|x-\xi\|^2),\varphi\in\mathcal{T},\xi\in\mathcal{D}(C)\}\subset\mathcal{D}(H_0).\tag{3.12}$$

Then $F$ approximates the metric $q(x,y)=\|x-y\|\wedge1$. In addition, for $\lambda>0$, if $f\in\mathcal{D}(H_0)$, then $\lambda f\in\mathcal{D}(H_0)$.

LEMMA 3.4.

(1) *For each* $f \in \mathcal{D}(H_0)$, *there exists* $f_n \in \mathcal{D}(H_n)$ *such that*

$$\sup_n \sup_x (|f_n(x)| + |H_n f_n(x)|) < \infty$$

*and*

$$\lim_{n\to+\infty} f_n(x_n) = f(x_0),$$
$$\limsup_{n\to+\infty} H_n f_n(x_n) \le (H_0 f)^*(x_0)$$

*whenever* $x_n \to x_0$.

(2) *For each* $f \in \mathcal{D}(H_1)$, *there exists* $f_n \in \mathcal{D}(H_n)$ *such that*

$$\sup_n \sup_x (|f_n(x)| + |H_n f_n(x)|) < \infty$$

*and*

$$\lim_{n\to+\infty} f_n(x_n) = f(x_0),$$
$$\liminf_{n\to+\infty} H_n f_n(x_n) \ge (H_1 f)_*(x_0)$$

*whenever* $x_n \to x_0$.

PROOF. Let us present the proof for $H_0$ only; the case for $H_1$ is similar. To further simplify, let us just verify the case for test functions in $\mathcal{D}(H_0)$ of the form

$$f(x) = g(x) + \varphi_2(h_{\varepsilon,y}(x)) \equiv \varphi_1(\|x-\xi\|^2) + \varphi_2(h_{\varepsilon,y}(x)) \in \mathcal{D}(H_0),$$

where $h_{\varepsilon,y}$ is defined by (3.8). By Condition 1.11, there exists $\xi_n \in E$ such that

$$\lim_{n\to\infty} \|\xi - \xi_n\| + \|C^0\xi - C_n\xi_n\| = 0.$$

Let $a_n > 1$ satisfy $\lim_{n\to\infty} a_n = +\infty$. We define $h_{n,\varepsilon,y}$ according to (3.9). Let

$$g_n(x) \equiv \varphi_1(\|x-\xi_n\|^2),$$
$$f_n(x) \equiv g_n(x) + \varphi_2(h_{n,\varepsilon,y}(x)) \in \mathcal{D}(H_n).$$

By part 3 of Lemma A.12,

$$\lim_{n\to\infty} \sup_{x\in K} |f_n(x) - f(x)| = 0 \qquad \text{for each } K \subset E \text{ compact.}$$

Apply (3.2), (A.18) and (A.20) to (3.6):

$$
\begin{aligned}
H_n f_n(x) \le\ & 2\varphi_1'(\|x-\xi_n\|^2)\langle x-\xi_n, C_n\xi_n\rangle \\
& + \sup_{r\ge 0}\varphi_2'(r) + \langle F_n(x), D(g_n+\varphi_2\circ h_{n,\varepsilon,y})(x)\rangle \\
& + \frac{1}{2}\|B_n^*(x)D(g_n+\varphi_2\circ h_{n,\varepsilon,y})(x)\|_{U_0}^2 \\
& + \frac{1}{2n}\operatorname{Tr}[D^2 g_n(x)B_n(x)B_n^*(x)] \\
& + \frac{1}{2n}\operatorname{Tr}[D^2(\varphi_2\circ h_{n,\varepsilon,y_1})(x)B_n(x)B_n^*(x)]
\end{aligned}
$$

$$
\begin{aligned}
\le\ & 2\varphi_1'(\|x-\xi_n\|^2)\langle x-\xi_n, C_n\xi_n\rangle + \sup_{r\ge 0}\varphi_2'(r) \\
& + \sup_{\|q\|\le\varphi_2'(h_{n,\varepsilon,y}(x))}\Big(\langle F_n(x), Dg_n(x)+q\rangle \\
& \qquad + \frac{1}{2}\|B_n^*(x)(Dg_n(x)+q)\|_{U_0}^2\Big) \\
& + \frac{1}{2n}\operatorname{Tr}[D^2 g_n(x)B_n(x)B_n^*(x)] \\
& + \frac{1}{2n}\operatorname{Tr}[D^2(\varphi_2\circ h_{n,\varepsilon,y})(x)B_n(x)B_n^*(x)],
\end{aligned} \tag{3.13}
$$

where

$$
Dg_n(x) = 2\varphi_1'(\|x-\xi_n\|^2)(x-\xi_n), \tag{3.14}
$$

$$
D^2 g_n(x) = 2\varphi_1'(\|x-\xi_n\|^2)I + 2\varphi_1''(\|x-\xi_n\|^2)(x-\xi_n)\otimes(x-\xi_n). \tag{3.15}
$$

Let $x_n \to x_0$, and denote

$$
\delta_n = \frac{1}{n}|||B_n(x_n)|||_{L_2(U_0,E)}^2.
$$

Then $\delta_n \to 0$ according to Condition 1.11(3). Taking $a_n = \delta_n^{-1/2}$ to be the one in (3.9), then $a_n \to +\infty$ and

$$
\frac{1}{n}a_n|||B_n(x_n)|||_{L_2(U_0,E)}^2 = a_n\delta_n = \delta_n^{1/2} \to 0.
$$

By (A.21),

$$
\lim_{n\to\infty}\frac{1}{2n}\operatorname{Tr}[D^2(\varphi_2\circ h_{n,\varepsilon,y})(x_n)B_n(x_n)B_n^*(x_n)] = 0. \tag{3.16}
$$

Therefore, by (3.13) through (3.16),

$$
\limsup_{n\to+\infty} H_n f_n(x_n) \le (H_0 f)^*(x_0)
$$

whenever $x_n \to x_0 \in E$. $\square$

3.4. *The comparison principle.* Let $\alpha > 0$, and let $h \in C_b(E)$ be uniformly continuous on $E$. The main goal of this subsection is to prove the comparison principle in Lemma 3.10.

In what follows, we extend the operator $H_0, H_1$ and connect Feng and Kurtz's definition of viscosity solution (Definition 1.14) with those in [31, 32] and [8]. We will introduce a new set of operators $\tilde{H}_0, \tilde{H}_1, \hat{H}_0, \hat{H}_1$ and will denote $\bar{H}_0$, $\bar{H}_1$ closures of $H_0, H_1$ under the graph norm topology in $B(E)$. We will clarify the relations among the next four sets of equations:

$$(I - \alpha H_0)f = h, \tag{3.17}$$

$$(I - \alpha H_1)f = h; \tag{3.18}$$

$$(I - \alpha \bar{H}_0)f = h, \tag{3.19}$$

$$(I - \alpha \bar{H}_1)f = h; \tag{3.20}$$

$$(I - \alpha \tilde{H}_0)f = h, \tag{3.21}$$

$$(I - \alpha \tilde{H}_1)f = h \tag{3.22}$$

and

$$(I - \alpha \hat{H}_0)f = h, \tag{3.23}$$

$$(I - \alpha \hat{H}_1)f = h. \tag{3.24}$$

Let $\varphi_i \in \mathcal{T}$ [see (1.32)], $y_i \in E$ and $\xi \in \mathcal{D}(C)$,

$$\begin{aligned} g(x) &= \varphi_1(\|x - \xi\|^2), \\ f(x) &= g(x) + \varphi_2(d_C(x, y_1)) + \cdots + \varphi_{k+1}(d_C(x, y_k)); \end{aligned}$$

we define single-valued operator

$$\begin{aligned} \tilde{H}_0 f(x) &= 2\varphi_1'(\|x - \xi\|^2)\langle x - \xi, C^0 \xi\rangle \\ &\quad + \Big(\sup_{r\ge 0} \varphi_2'(r) + \cdots + \sup_{r \ge 0} \varphi_{k+1}'(r)\Big) \\ &\quad + \sup_{\|q\| \le \varphi_2'(d_C(x,y_1)) + \cdots + \varphi_{k+1}'(d_C(x,y_k))} (\langle F(x), Dg(x) + q\rangle \\ &\qquad + \tfrac{1}{2}\|B^*(x)(Dg(x) + q)\|_{U_0}^2). \end{aligned} \tag{3.25}$$

By item (1) in Lemma A.12, (3.10) is equal to (3.25) when $\|x\|$ is sufficiently large, independent of the $\varepsilon$. In addition, by (A.17),

$$\lim_{\varepsilon \to 0+} \sup_{x \in E} |h_{\varepsilon, y}(x) - d_C(x, y)| = 0.$$

Therefore sending $\varepsilon \to 0$, we obtain $\tilde{H}_0 \subset \bar{H}_0$, where $\bar{H}_0$ is the closure of $H_0$ under the uniform norm for $B(E)$. Similarly, let

$$g(x) = -\varphi_1(\|x-\xi\|^2),$$
$$f(x) = g(x) - (\varphi_2(d_C(x,y_1)) + \cdots + \varphi_{k+1}(d_C(x,y_k)))$$

and define

$$\begin{aligned}\tilde{H}_1 f(x) = {} & -2\varphi_1'(\|x-\xi\|^2)\langle x-\xi, C^0\xi\rangle \\ & - \Big(\sup_{r\ge 0}\varphi_2'(r) + \cdots + \sup_{r\ge 0}\varphi_{k+1}'(r)\Big) \\ & + \inf_{\|q\|\le \varphi_2'(d_C(x,y_1))+\cdots+\varphi_{k+1}'(d_C(x,y_k))} (\langle F(x), Dg(x)+q\rangle \\ & \qquad + \tfrac{1}{2}\|B^*(x)(Dg(x)+q)\|_{U_0}^2).\end{aligned} \tag{3.26}$$

Then $\tilde{H}_1 \subset \bar{H}_1$.

LEMMA 3.5. *$\overline{f}$ is a viscosity subsolution of (3.17) for $H_0$ if and only if it is a viscosity subsolution of (3.19) for $\bar{H}_0$; both imply $\overline{f}$ is also a viscosity subsolution of (3.21) for $\tilde{H}_0$.*

*$\underline{f}$ is a viscosity supersolution of (3.18) for $H_1$ if and only if it is a viscosity supersolution of (3.20) for $\bar{H}_1$; both imply $\underline{f}$ is also a viscosity supersolution of (3.22) for $\tilde{H}_1$.*

*Hence, the comparison principle for subsolution of (3.21) and supersolution of (3.22) implies the comparison principles for (3.17) and (3.18), as well as those for (3.19) and (3.20).*

*In this lemma viscosity solution is always meant in the sense of Definition 1.14.*

PROOF. The conclusion follows from the fact that $\tilde{H}_0 \subset \bar{H}_0$, $\tilde{H}_1 \subset \bar{H}_1$; and the definition of viscosity solution in Definition 1.14. □

We discuss some properties enjoyed by functions in $\mathcal{D}(\tilde{H}_i)$, $i = 0, 1$.

LEMMA 3.6. *Let $f_0 \in \mathcal{D}(\tilde{H}_0)$. Suppose $x_0 \in E$ satisfies $(\overline{f} - f_0)(x_0) = \sup_{x\in E}(\overline{f} - f_0)(x)$. Let $0 \le \varphi \in C^2([0,\infty))$ be nondecreasing, $\varphi(r) = r$ when $r \le 1$ and $\varphi(r) = 2$ when $r \ge 2$. Let $\theta > 0$. We introduce perturbation of $f_0$:*

$$f_\theta(x) = f_0(x) + \theta\varphi(d_C(x,x_0)). \tag{3.27}$$

*Then $f_\theta$ has the following properties:*

(a) *$f_\theta \in \mathcal{D}(\tilde{H}_0)$ and*

$$(\overline{f} - f_\theta)(x_0) > (\overline{f} - f_\theta)(x), \qquad x \ne x_0.$$

(b) *For any* $\{x_n\} \subset E$ *satisfying*

$$\lim_{n\to\infty}(\overline{f} - f_\theta)(x_n) = \sup_{x\in E}(\overline{f} - f_\theta)(x),$$

*we have* $x_n \to x_0$ *and* $\overline{f}(x_n) \to \overline{f}(x_0)$.

(c)

$$\limsup_{\theta\to 0+}(\tilde{H}_0 f_\theta)^*(x_0) \le (\tilde{H}_0 f_0)^*(x_0),$$

$$\liminf_{\theta\to 0+}(\tilde{H}_1 f_\theta)_*(x_0) \ge (\tilde{H}_1 f_0)_*(x_0).$$

PROOF. Part (a) follows from the definition of $f_\theta$.

We prove part (b) next. By (a),

$$\lim_n(\overline{f} - f_\theta)(x_n) = \sup_x(\overline{f} - f_\theta)(x) = (\overline{f} - f_\theta)(x_0) = (\overline{f} - f_0)(x_0).$$

Therefore

$$\begin{aligned}(\overline{f} - f_0)(x_0) = (\overline{f} - f_\theta)(x_0) &= \lim_n(\overline{f} - f_\theta)(x_n)\\ &= \lim_n(\overline{f} - f_0)(x_n) - \theta\varphi(d_C(x_n, x_0))\\ &\le \liminf_n(\overline{f} - f_0)(x_n) = (\overline{f} - f_0)(x_0).\end{aligned}$$

Hence

$$\begin{aligned}\lim_{n\to\infty}\theta d_C(x_n, x_0) &= \lim_{n\to\infty}\{(\overline{f} - f_0)(x_n) - ((\overline{f} - f_0)(x_n) - \theta\varphi(d_C(x_n, x_0)))\}\\ &= 0,\end{aligned}$$

which implies $x_n \to x_0$ and $(\overline{f} - f_0)(x_n) \to (\overline{f} - f_0)(x_0)$.

Part (c) follows from direct verification. □

LEMMA 3.7.

(a) *If* $\overline{f}$ *is a viscosity subsolution of* (3.21) *in the sense of Definition* 1.14*, then it is also a viscosity subsolution in the sense of Definition* 1.16 *for* $\tilde{H}_0$.

(b) *If* $\underline{f}$ *is a viscosity supersolution of* (3.22) *in the sense of Definition* 1.14*, then it is also a viscosity supersolution in the sense of Definition* 1.16 *for* $\tilde{H}_1$.

PROOF. We prove part (a) only. The proof for part (b) is similar. Let $\varphi_i \in \mathcal{T}$ [see (1.32)], $y_i \in E$ and $\xi \in \mathcal{D}(C)$. Consider

$$g_0(x) = \varphi_1(\|x - \xi\|^2),$$

$$f_0(x) = g_0(x) + \varphi_2(d_C(x, y_1)) + \cdots + \varphi_{k+1}(d_C(x, y_k)) \in \mathcal{D}(\tilde{H}_0),$$

and $x_0 \in E$ such that $(\overline{f} - f_0)(x_0) = \sup_{x\in E}(\overline{f} - f_0)(x)$. Define $f_\theta$ according to (3.27). By Lemma 3.6, $f_\theta \in \mathcal{D}(\tilde{H}_0)$ and

$$(\overline{f} - f_\theta)(x_0) > (\overline{f} - f_\theta)(x), \qquad x \neq x_0.$$

Since $\overline{f}$ is a subsolution of (3.21) in the sense of Definition 1.14, there exists a sequence $\{x_n\} \subset E$ such that

$$\lim_n (\overline{f} - f_\theta)(x_n) = \sup_x (\overline{f} - f_\theta)(x)$$

and

$$\limsup_{n\to\infty} (\alpha^{-1}(\overline{f} - h)(x_n) - (\tilde{H}_0 f_\theta)^*(x_n)) \leq 0.$$

By Lemma 3.6, $x_n \to x_0$, $\overline{f}(x_n) \to \overline{f}(x_0)$ and

$$\limsup_{\theta\to 0+} (\tilde{H}_0 f_\theta)^*(x_0) \leq (\tilde{H}_0 f_0)^*(x_0).$$

Hence

$$\alpha^{-1}(\overline{f} - h)(x_0) \leq \limsup_{\theta\to 0+} (\tilde{H}_0 f_\theta)^*(x_0) \leq (\tilde{H}_0 f_0)^*(x_0).$$

□

We define $\hat{H}_0$ and $\hat{H}_1$ below. For

$$\begin{aligned} g(x) &= \frac{\mu}{2}\|x - \xi\|^2 \qquad && \forall \mu > 0, \xi \in \mathcal{D}(C), \\ f(x) &= g(x) + \rho d_C(x, y) && \forall y \in E, \rho > 0 \end{aligned} \tag{3.28}$$

(recall the definition of $d_C$ in Definition 1.13), we define

$$\begin{aligned} \hat{H}_0 f(x) &= \mu\langle x - \xi, C^0\xi\rangle + \rho \\ &\quad + \sup_{\|q\|\leq\rho} (\langle F(x), Dg(x) + q\rangle + \tfrac{1}{2}\|B^*(x)(Dg(x) + q)\|^2_{U_0}). \end{aligned} \tag{3.29}$$

Similarly, for

$$\begin{aligned} g(x) &= -\frac{\mu}{2}\|x - \xi\|^2 \qquad && \forall \mu > 0, \xi \in \mathcal{D}(C), \\ f(x) &= g(x) - \rho d_C(x, y) && \forall y \in E, \rho > 0, \end{aligned} \tag{3.30}$$

we define

$$\begin{aligned} \hat{H}_1 f(x) &= -\mu\langle x - \xi, C^0\xi\rangle - \rho \\ &\quad + \inf_{\|q\|\leq\rho} (\langle F(x), Dg(x) + q\rangle + \tfrac{1}{2}\|B^*(x)(Dg(x) + q)\|^2_{U_0}). \end{aligned} \tag{3.31}$$

$\tilde{H}_0, \tilde{H}_1$ are local operators, therefore we can get rid of the localization functions $\varphi_k \in \mathcal{T}$ to arrive at $\hat{H}_0, \hat{H}_1$. We omit the proof here. Such argument is standard. In Lemma 6.1 and Theorem 6.1 of [20], these types of arguments are used to prove equivalence of different definitions of viscosity solution. We have the following conclusion.

LEMMA 3.8. *If $\overline{f}$ is the viscosity subsolution of (3.21) for $\tilde{H}_0$, then it is also the subsolution of (3.23) for $\hat{H}_0$; both in the sense of Definition 1.16.*

*If $\underline{f}$ is the viscosity supersolution of (3.22) for $\tilde{H}_1$, then it is also the supersolution of (3.24) for $\hat{H}_1$; both in the sense of Definition 1.16.*

Summarizing conclusions in Lemmas 3.5, 3.7 and 3.8, we have the next result.

LEMMA 3.9. *Let $\overline{f}$ be a subsolution to (3.17) for $H_0$ in the sense of Definition 1.14 (Feng and Kurtz); then it is a subsolution to (3.23) for $\hat{H}_0$ in the sense of Definition 1.16 (Tataru–Crandall–Lions).*

*Let $\underline{f}$ be a supersolution to (3.18) for $H_1$ in the sense of Definition 1.14 (Feng and Kurtz); then it is a supersolution to (3.24) for $\hat{H}_1$ in the sense of Definition 1.16 (Tataru–Crandall–Lions).*

We will study the comparison principle for viscosity solutions in the sense of Definition 1.16 in Section 5. In view of Lemma 3.9, Theorem 5.1 implies the following.

LEMMA 3.10. *Let $\overline{f}$ be a subsolution to (3.17) for $H_0$ and let $\underline{f}$ be a supersolution to (3.18) for $H_1$, both in the sense of Definition 1.14. Then $\overline{f} \le \underline{f}$.*

3.5. *The large deviation theorem.*

THEOREM 3.11. *Suppose Conditions 1.11 and 1.12 are satisfied. Let $X_n \in C_E[0,\infty)$ be the solution of (1.7). Suppose further that $\{X_n(0)\}$ satisfies the large deviation principle with good rate function $I_0$ on $E$.*

*Then:*

(a) *$\{X_n\}$ is exponentially tight;*
(b) *the following limit exists and defines an operator semigroup on $C_b(E)$:*

$$V(t)f(x) = \lim_{n\to\infty} \frac{1}{n} \log E[e^{nf(X_n(t))}|X_n(0) = x]; \tag{3.32}$$

(c) *the large deviation principle holds for $\{X_n\}$ with good rate function $I$:*

$$I(x) = I_0(x(0)) + \sup_{0\le t_1\le\cdots\le t_m} \left(\sum_{i=1}^m I_{t_i - t_{i-1}}(x(t_i)|x(t_{i-1}))\right), \tag{3.33}$$

*where*

$$I_t(y|x) = \sup_{f\in C_b(E)} (f(y) - V(t)f(x)).$$

PROOF. Define $F \subset C_b(E)$ according to (3.12). The operator convergence in Lemma 3.4 and the comparison principle in Lemma 3.10 imply that Theorem 2.3 holds. Consequently the conclusion follows. $\square$

**4. Application to special cases.** We solve Examples 1.2, 1.5 and 1.8 as special cases of Theorem 3.11.

4.1. *Stochastic Allen–Cahn equation.* Recall that we take $E = U_0 = L^2(\mathcal{O})$. Let $\omega = \sup_r |V''(r)|$. We take

$$(Cx)(\theta) \equiv \Delta x(\theta) - V'(x(\theta)) + V'(0) - \omega x(\theta) \tag{4.1}$$

where

$$\mathcal{D}(C) = H^2(\mathcal{O}) \equiv \left\{ x : x, \frac{\partial}{\partial \theta_i} x, \frac{\partial^2}{\partial \theta_i \, \partial \theta_j} x \in L^2(\mathcal{O}), i, j = 1, \dots, d \right\}$$

and

$$F(x) = -V'(0) + \omega x.$$

$C$ is $m$-dissipative by the usual theory of semilinear equation. Recall the $C_n$ and $F_n$ in (1.18); we have

$$\begin{aligned} &\lim_{n\to\infty} \|C_n \xi - C\xi\|_{L^2(\mathcal{O})} = 0 \qquad \forall \xi \in \mathcal{D}(C), \\ &\lim_{x_n \to x_0} F_n(x_n) = F(x_0). \end{aligned} \tag{4.2}$$

In addition, let $B_n(x)$ be defined according to (1.16); then

$$\lim_{x_n \to x_0, p_n \to p_0} \|B_n^*(x_n) p_n\|_{L^2(\mathcal{O})} = \|B^*(x_0) p_0\|_{L^2(\mathcal{O})}.$$

Let $\{e_1, \dots, e_k, \dots\}$ be the orthonormal system for $U_0 = L^2(\mathcal{O})$ as defined in (1.11). Let $\sigma, \varphi$ be defined according to (1.9). Since

$$\begin{aligned} \|\!|B_n(x_n)|\!\|^2_{L_2(U_0,E)} &= \mathrm{Tr}(B_n^*(x_n) B_n(x_n)) = \sum_k \|B_n(x_n) e_k\|^2 \\ &= \sum_k \sum_{i=(1,\dots,1)}^{(m_n,\dots,m_n)} \langle \varphi(\cdot, \langle P_n x_n, \xi\rangle) e_k, e_i\rangle^2 \\ &= \sum_{i=(1,\dots,1)}^{(m_n,\dots,m_n)} \sum_k \langle \varphi(\cdot, \langle P_n x_n, \xi\rangle) e_i, e_k\rangle^2 \\ &= \sum_{i=(1,\dots,1)}^{(m_n,\dots,m_n)} \|\varphi(\cdot, \langle P_n x_n, \xi\rangle) e_i\|^2 \le m_n^d \sup_{\theta, r} \varphi^2(\theta, r), \end{aligned}$$

Condition 1.11(3) holds under the scaling requirement (1.14).

Finally, Condition 1.12 is verified by Lemma A.5. Therefore, Theorem 1.4 follows from Theorem 3.11.

To simplify the form of the rate function from (1.31) to a time integral form as in (1.19), we need additional work. The basic idea is the same as that presented in Example 2.5, with some technical complications because of the infinite-dimensional state space. A general result for rate function representation is developed in Chapter 8 of [18]. Applying such result, representation (1.19) for lattice versions of the stochastic Allen–Cahn equation is proved rigorously in Chapter 13 of [18]. This procedure can be carried out similarly here. Below, we only provide a sketch.

First, the form of $\hat{H}_0$, $\hat{H}_1$ in (3.29) and (3.31) induces an optimal controlled PDE problem:

$$\begin{aligned}\frac{\partial}{\partial t}x(t,\theta) &= Cx + F(x) + B(x)u(t)\\ &= \Delta x(t,\theta) - V'(x(t,\theta)) + \sigma(x,\theta)u(t,\theta)\end{aligned} \tag{4.3}$$

which is well defined under a finite *running cost* assumption:

$$\tfrac{1}{2}\int_0^\infty \int_{\mathcal{O}} u^2(t,\theta)\,d\theta\,dt < \infty.$$

By the dynamic programming principle and the comparison principle for (3.23) and (3.24), we can prove that the $V(t)$ in (3.32) has the form

$$\begin{aligned}V(t)f(x_0) = \sup\bigg\{ & f(x(t))\\ & - \tfrac{1}{2}\int_0^t \int_{\mathcal{O}} u^2(s,\theta)\,d\theta\,ds : (x,u) \text{ satisfies (4.3) and } x(0) = x_0\bigg\}.\end{aligned}$$

Then from this, we derive

$$\begin{aligned}I_t(x_1|x_0) = \inf\bigg\{ & \int_0^t \int_{\mathcal{O}} \tfrac{1}{2}u^2(t,\theta)\,d\theta\,ds|\\ & (x,u) \text{ satisfies (4.3) with } x(0) = x_0, x(t) = x_1\bigg\}.\end{aligned}$$

Combine the above form with (4.3); the variational form of $I(x)$ in (1.19) follows.

4.2. *Stochastic Cahn–Hilliard equation.* Let $\omega \equiv \sup_r |V''(r)|^2/4$. We define

$$(Cx)(\theta) \equiv \Delta(-\Delta x(\theta) + V'(x(\theta))) - \omega x(\theta) \tag{4.4}$$

for

$$x \in \mathcal{D}(C) = W^{4,2}(\mathcal{O})$$
$$\equiv \left\{ x : x, \frac{\partial}{\partial \theta_i} x, \frac{\partial^2}{\partial \theta_i \, \partial \theta_j} x, \frac{\partial^3}{\partial \theta_i \, \partial \theta_j \, \partial \theta_k} x, \frac{\partial^4}{\partial \theta_i \, \partial \theta_j \, \partial \theta_k \, \partial \theta_l} x \in L^2(\mathcal{O}) \right\}$$

and

$$(F(x))(\theta) \equiv \omega x(\theta), \qquad x \in L^2(\mathcal{O}).$$

We recall that the $C_n, F_n$ are defined as in (1.23). By Lemma A.1, $C$ and $C_n$ are $m$-dissipative in $L^2(\mathcal{O})$. Furthermore, the type of convergence in (4.2) holds here by direct verification.

By Lemma A.7, Condition 1.12 is also satisfied.

Theorem 1.7 follows from Theorem 3.11.

The rate function representation in (1.24) can be proved using similar arguments as in the Allen–Cahn case. The controlled PDE becomes

$$\frac{\partial}{\partial t} x(t,\theta) = Cx + F(x) + B(x)u(t) = \Delta(-\Delta x(t,\theta) + V'(x(t,\theta))) + u(t,\theta).$$

The running cost structure is the same.

4.3. *Stochastic quasilinear equation with viscosity.* Let $C_n$ be defined according to (1.27) and let

$$Cx = \alpha \Delta_\theta x - \partial_\theta \phi(x) - \omega x, \qquad x \in H^2(\mathcal{O}). \tag{4.5}$$

Then both $C_n$ and $C$ are $m$-dissipative operators in $L^2(\mathcal{O})$ (Lemma A.4).

The compact containment estimate is provided in Lemma A.9. Following the same arguments as above, Theorem 1.9 follows as a special case of Theorem 3.11.

Rate function representation is the same as the Allen–Cahn case. The controlled PDE is

$$\frac{\partial}{\partial t} x(t,\theta) = Cx + F(x) + B(x)u(t) = \alpha \, \partial^2_{\theta\theta} x(t,\theta) - \partial_\theta \phi(x(t,\theta)) + u(t,\theta).$$

**5. A class of Hamilton–Jacobi equation in Hilbert space.** The purpose of this section is to present a self-contained proof of the comparison principle for (5.2) and (5.3) in the viscosity solution sense by Tataru–Crandall–Lions (Definition 1.16)—Theorems 5.1 and 5.2. The whole section is independent of the rest of the paper and can be read separately.

We point out that the comparison results in [31, 32] and in [8] cannot be directly borrowed here, because the $\hat{H}_0, \hat{H}_1$ are not exactly of the same form as considered there. For example, when defining these operators, we restrict the domains and use rougher estimates than the $D_C^+$ and $D_C^-$ (Definition 3.3).

This allows us to relate $\hat{H}_0, \hat{H}_1$ with other operators which arise as the kind of limits required by Theorem 2.3 with graphs contained in $C_b(E) \times B(E)$. Second but more importantly, the quadratic nonlinearity in (5.6) is worse than that assumed in (ii) of (49) in [8]. We explore convexity to cure this problem. Despite these differences, the main ideas of [8, 31, 32], still apply and all we need are modifications and refinements at various places. We follow [8] and present the proof through a doubling technique—Lemma 5.9. To make this paper self-contained, we will repeat the important steps of [8] and omit minor details. We make detailed references for all omitted steps so that they can easily be recovered if needed.

Let $E$ be a separable real Hilbert space, let $C \subset E \times E$ be a possibly multivalued nonlinear $m$-dissipative operator with $\overline{\mathcal{D}(C)} = E$ and $(0,0) \in C$. Further suppose function $F : E \to E$ and operator $B(x) : E \to E$ for each $x \in E$ satisfy

$$L_{F,B} \equiv \sup_{x \neq y} \frac{\|F(x) - F(y)\| + \|\|B(x) - B(y)\|\|}{\|x - y\|} < \infty. \tag{5.1}$$

Recall the definition of $\hat{H}_0, \hat{H}_1 \subset C(E) \times M(E)$ in (3.29) and (3.31); we consider the following Hamilton–Jacobi equations written in the resolvent form: let $h \in C_b(E)$ and $\alpha > 0$,

$$(I - \alpha \hat{H}_0) f = h \tag{5.2}$$

and

$$(I - \alpha \hat{H}_1) f = h. \tag{5.3}$$

We prove the following.

Theorem 5.1 (Comparison principle). *Let $\overline{f}$ be a subsolution of (5.2) and $\underline{f}$ be a supersolution of (5.3), both in the sense of Definition 1.16.*

*Suppose $h$ is uniformly continuous and (5.1) is satisfied. Then*

$$\overline{f} \leq \underline{f}.$$

Indeed, we will prove a theorem covering more general situations. Let $G(x, p) \in C(E \times E)$. We define single-valued operators $\hat{H}_0, \hat{H}_1 \subset C(E) \times M(E)$: for each $f$ in (3.28), we define

$$\hat{H}_0 f(x) \equiv \mu \langle x - \xi, C^0 \xi \rangle + \rho + \sup_{\|q\| \leq \rho} G(x, Dg(x) + q). \tag{5.4}$$

For each $f$ in (3.30), we define

$$\hat{H}_1 f(x) \equiv -\mu \langle x - \xi, C^0 \xi \rangle - \rho + \inf_{\|q\| \leq \rho} G(x, Dg(x) + q). \tag{5.5}$$

The operators in (5.2) and (5.3) are special cases of the above ones with $G$ given by

$$G(x,p) = \langle F(x), p\rangle + \tfrac{1}{2}\|B^*(x)p\|^2_{U_0}. \tag{5.6}$$

THEOREM 5.2. *Suppose $\alpha > 0$, $h$ is uniformly continuous and Condition 5.5 is satisfied for $G$. Define $\hat{H}_0, \hat{H}_1$ according to (5.4) and (5.5).*

*Let $\overline{f}$ be a subsolution of (5.2) and $\underline{f}$ be a supersolution of (5.3), both in the sense of Definition 1.16.*

*Then*

$$\overline{f} \le \underline{f}.$$

5.1. *Perturbed optimization principle.* There is no a priori guarantee that the extrema in (1.40) and (1.41) of Definition 1.16 always exist. We have to carefully choose test functions $f_0$ to make sure that the definition is not an empty one. Ekeland's perturbed optimization principle [14] claims that, if we add a small perturbation to the test function using the norm of the Hilbert space, we can always attain the extrema. If we apply this technique in the viscosity solution context, we also want the perturbed $Hf_0$ to be close to the unperturbed one, so that the equation we consider does not change much. These considerations lead to the Tataru distance function $d_C$ in Definition 1.13.

The following is adapted from Proposition 2.1 of [31], which generalizes Ekeland's principle. In the adaptation, we have taken $g = -u$, where $u$ is the function in the original proposition.

LEMMA 5.3. *Let $K$ be an abstract set and $B : K \times K \to [0, +\infty)$, with the following properties:*

(a) *$B(x,x) = 0$ for all $x \in K$;*
(b) *$B(x,y) + B(y,z) \ge B(x,z)$ for all $x,y,z \in K$;*
(c) *for each $\{x_n\} \subset K$ satisfying $\sum_{n=1}^{\infty} B(x_n, x_{n+1}) < \infty$, there exists $x \in K$ such that $\lim_{n\to+\infty} B(x_n, x) = 0$.*

*Let $g : K \to [-\infty, +\infty)$, $\sup_{x\in K} g(x) < +\infty$. Furthermore, if $\{x_n\} \subset K$ and $x \in K$ satisfy $\sum_{n=1}^{\infty} B(x_n, x_{n+1}) < +\infty$ and $\lim_{n\to+\infty} B(x_n, x) = 0$, then $g(x) \ge \limsup g(x_n)$.*

*Then, for each $\varepsilon > 0$, and $x_0$ such that $g(x_0) \ne -\infty$, there exists $x_\varepsilon$ such that:*

(1) $g(x_0) + \varepsilon B(x_0, x_\varepsilon) \le g(x_\varepsilon)$,
(2) $g(x) - \varepsilon B(x, x_\varepsilon) \le g(x_\varepsilon)$, $x \in K$.

Remark 5.4. Let function $g : E \to R$ be bounded and upper semicontinuous in the norm topology of $E$. If we take $K = E$ and $B = d_C$, then the assumptions regarding $B$ are satisfied (Proposition 2.2 of [31]), and:

(a) For each $x_0 \in E$ and $\varepsilon > 0$, there exists an $x_1 \in E$ such that

$$g(x_0) + \varepsilon d_C(x_0, x_1) \le g(x_1) \tag{5.7}$$

and

$$g(x) - \varepsilon d_C(x, x_1) \le g(x_1), \qquad x \in E. \tag{5.8}$$

(b) Let $\varepsilon > 0$ and $x_0 \in E$ be such that

$$\sup_{x\in E} g(x) \le g(x_0) + \varepsilon^2.$$

Then there exists $x_1 \in E$ such that not only (5.7) and (5.8) hold, but also

$$d_C(x_0, x_1) \le \varepsilon.$$

Part (a) is a consequence of the above proposition. See also Lemmas 2.4, 2.5 of [8]. Part (b) follows from (5.7):

$$\varepsilon d_C(x_0, x_1) \le g(x_1) - g(x_0) \le g(x_1) - \sup_{x\in E} g(x) + \varepsilon^2 \le \varepsilon^2.$$

Similarly, an analogous result holds when we take $K = E \times E$ and $B((x_1, y_1), (x_2, y_2)) = d_C(x_1, x_2) + d_C(y_1, y_2)$. We will need such result for (5.14).

5.2. *The comparison principle.* We make the following structural assumption about $G$.

Condition 5.5.

(1) $G(x, p) \in C(E \times E)$; for each $\lambda > 1$ and $M > 0$ fixed, there exist $\rho_\lambda(r)$, $\sigma_M(r) \in C(R^+)$ with $\rho_\lambda(0) = 0$ and $\sigma_M(0) = 0$ such that

$$\lambda G\left(x, \frac{\mu(x-y)}{\lambda}\right) - G(y, \mu(x-y)) \le \rho_\lambda(\|x-y\| + \mu\|x-y\|^2), \qquad \mu > 0,$$

and

$$\sup_{\|x\|+\|p\|+\|q\|<M, \|q-p\|\le r} |G(x, p) - G(x, q)| \le \sigma_M(r).$$

(2) There exists a nondecreasing function $0 \le \varphi \in C^1([0,\infty))$ slowly growing to infinity in the sense that $\lim_{r\to\infty} \varphi(r) = \infty$ and $\sup_{r\ge 0} |r\varphi'(r^2)| + |r\varphi'(r)| < \infty$. $\varphi' > 0$.

For $M > 0$, $\lambda > 1$, there exists $\gamma_{\lambda,M}(r) \in C(R^+)$ with $\gamma_{\lambda,M}(0) = 0$ such that

$$\sup_{x\in E, \|p\|\le M} \lambda G\left(x, \frac{p + \varepsilon\varphi'(\|x\|^2)x}{\lambda}\right) - G(x, p) \le \gamma_{\lambda,M}(\varepsilon).$$

REMARK 5.6. Condition 5.5 implies that for any $1<\lambda_0<\lambda_1$,

$$\lambda_1 G\Big(x, \frac{\mu(x-y)}{\lambda_1}\Big) - \lambda_0 G\Big(y, \frac{\mu(x-y)}{\lambda_0}\Big) \\ \le \lambda_0 \rho_{\lambda_1/\lambda_0}\Big(\|x-y\| + \frac{\mu}{\lambda_0}\|x-y\|^2\Big)$$

and

$$\sup_{x\in E, \|p\|\le M}\Big(\lambda_1 G\Big(x, \frac{p+\varepsilon\varphi'(\|x\|^2)x}{\lambda_1}\Big) - \lambda_0 G\Big(x, \frac{p}{\lambda_0}\Big)\Big) \le \lambda_0 \gamma_{\lambda_1/\lambda_0, M}\Big(\frac{\varepsilon}{\lambda_0}\Big).$$

$\varphi(r)=\log(1+r)$ satisfies the growth requirement in the condition. For many examples, such choice is good enough.

LEMMA 5.7. *Suppose (5.1) is satisfied; then the $G$ in (5.6) satisfies Condition 5.5.*

PROOF. To verify this, let $\lambda>1$; noting

$$\Big(\frac{1}{\lambda}-1\Big)x^2 + \frac{2}{\lambda}bx + \frac{1}{\lambda}b^2 \le \frac{b^2}{\lambda-1},$$

we have

$$\begin{aligned}\lambda G\Big(x, \frac{p}{\lambda}\Big) - G(y,p) &= \frac{1}{2}\Big(\frac{1}{\lambda}\|B^*(x)p\|_{U_0}^2 - \|B^*(y)p\|_{U_0}^2\Big) + \langle F(x)-F(y), p\rangle \\ &\le \frac{1}{2}\Big(\frac{1}{\lambda}\|(B^*(x)-B^*(y))p\|_{U_0}^2 + \Big(\frac{1}{\lambda}-1\Big)\|B^*(y)p\|_{U_0}^2 \\ &\qquad + \frac{2}{\lambda}\|(B^*(x)-B^*(y))p\|_{U_0}\|B^*(y)p\|_{U_0}\Big) \\ &\quad + L_{F,B}\|x-y\|\|p\| \\ &\le \frac{1}{2}\frac{\|(B^*(x)-B^*(y))p\|_{U_0}^2}{\lambda-1} + L_{F,B}\|x-y\|\|p\| \\ &\le \frac{1}{2}\frac{(L_{F,B}\|x-y\|\|p\|)^2}{\lambda-1} + L_{F,B}\|x-y\|\|p\|,\end{aligned}$$

where $L_{F,B}$ is the one in (5.1). We can take

$$\rho_\lambda(r) = \frac{1}{2}\frac{L_{F,B}^2}{\lambda-1}r^2 + L_{F,B}r.$$

Similarly, denoting $C_0 = \|\|B(0)\|\|$,

$$\lambda G\Big(x, \frac{p+\varepsilon q}{\lambda}\Big) - G(x,p) \le \frac{1}{2}\Big(\Big(\frac{1}{\lambda}-1\Big)\|B^*(x)p\|_{U_0}^2 + \frac{1}{\lambda}\|B^*(x)\varepsilon q\|_{U_0}^2$$

$$+\frac{2}{\lambda}\|B^*(x)p\|_{U_0}\|B^*(x)\varepsilon q\|_{U_0}\Big)$$

$$+\varepsilon\|F(x)\|\|q\|$$

$$\le \frac{1}{2}\frac{\|B^*(x)\varepsilon q\|^2_{U_0}}{\lambda-1}+\varepsilon L_{F,B}\|x\|\|q\|$$

$$\le \frac{\varepsilon^2}{2}\frac{(C_0\|q\|+L\|q\|\|x\|)^2}{\lambda-1}+\varepsilon L_{F,B}\|x\|\|q\|.$$

We can take $\varphi(r)=\log(1+r)$ and

$$\gamma_{\lambda,M}(\varepsilon)=\frac{\varepsilon^2}{2}\frac{(C_0+L)^2}{\lambda-1}+\varepsilon L_{F,B}.$$

Finally, for each $M>0$, we denote

$$C_1(M)=\sup_{\|x\|<M}|\!|\!|B(x)|\!|\!|,$$

$$C_2(M)=\sup_{\|x\|<M}|\!|\!|B(x)|\!|\!|^2M,$$

$$C_3(M)=\sup_{\|x\|<M}\|F(x)\|.$$

Then for $\|x\|+\|p\|+\|q\|<M$,

$$|G(x,p)-G(x,q)|\le (\|B^*(x)(p-q)\|^2_{U_0}$$
$$+2\|B^*(x)q\|_{U_0}\|B^*(x)(p-q)\|_{U_0}+\|F(x)\|\|p-q\|)$$
$$\le C_1^2(M)\|p-q\|^2+2C_2(M)\|p-q\|+C_3(M)\|p-q\|.$$

We can take $\sigma_M(r)=C_1^2(M)r^2+(2C_2(M)+C_3(M))r$. □

Now, we prove Theorem 5.2 in several steps.

We endow the product space $E\times E$ with inner product

$$\langle (x,y),(\xi,\eta)\rangle=\langle x,\xi\rangle+\langle y,\eta\rangle.$$

Denote

$$\mathcal{C}(x,y)=(Cx,Cy).$$

By the $m$-dissipativity of $C$, $\mathcal{C}$ induces a semigroup:

$$\mathcal{S}(t)(x,y)=(S(t)x,S(t)y),\qquad (x,y)\in E\times E.$$

For $\lambda>0$, $(x,y),(p,q)\in E\times E$, we define

$$G_\lambda((x,y);(p,q))\equiv \lambda G\Big(x,\frac{p}{\lambda}\Big)-G(y,-q),$$

and define a single-valued operator $H_{2,\lambda} \subset C(E \times E) \times M(E \times E)$ next. Let $\mathcal{D} = \mathcal{D}(H_{2,\lambda})$ consist of functions $\Phi$ defined as follows:

$$\phi(x,y) = \frac{\mu}{2}\|x-y\|^2 + \frac{\gamma}{2}(\varphi(\|x\|^2) + \varphi(\|y\|^2)), \qquad \mu, \gamma > 0, \tag{5.9}$$

where $\varphi \in C^1([0,\infty))$, $\varphi' \geq 0$, $\lim_{r\to+\infty} \varphi(r) = +\infty$; and

$$\Phi(x,y) = \phi(x,y) + \rho(d_C(x,x_0) + d_C(y,y_0)), \qquad \rho > 0, x_0, y_0 \in E.$$

We define, for each $\Phi$,

$$\begin{aligned} H_{2,\lambda}\Phi(x,y) &= D_{\mathcal{C}}^+\Phi(x,y) \\ &\quad + \sup_{\|p\|^2+\|q\|^2 \leq 2\rho^2} G_\lambda((x,y); D\phi(x,y) + (p,q)). \end{aligned} \tag{5.10}$$

We note that $H_{2,\lambda}\Phi(x,y)$ may be an unbounded function on $E \times E$.

We have a perturbation result.

Lemma 5.8. *Let $u \in B(E \times E)$ be upper semicontinuous, $\Phi_0(x,y) \in \mathcal{D}$ and $(\hat{x},\hat{y}) \in E \times E$ satisfy $(u - \Phi_0)(\hat{x},\hat{y}) = \sup_{(x,y)\in E\times E}(u-\Phi_0)(x,y)$. Suppose $\kappa > 0$; we define*

$$\Phi_\kappa(x,y) = \Phi_0(x,y) + \kappa(d_C(x,\hat{x}) + d_C(y,\hat{y})). \tag{5.11}$$

*Then $\Phi_\kappa$ has the following properties:*

(a)

$$(u - \Phi_\kappa)(\hat{x},\hat{y}) > (u-\Phi_\kappa)(x,y), \qquad (x,y) \neq (\hat{x},\hat{y}).$$

(b) *For any $\{(x_n,y_n)\} \subset E \times E$ satisfying*

$$\lim_{n\to\infty}(u-\Phi_\kappa)(x_n,y_n) = \sup_{(x,y)\in E\times E}(u-\Phi_\kappa)(x,y),$$

*we have $(x_n,y_n) \to (\hat{x},\hat{y})$ and $u(x_n,y_n) \to u(\hat{x},\hat{y})$.*

Proof. The same arguments as in the proof of Lemma 3.6 apply here. □

Lemma 5.9 (Doubling lemma). *Let $\lambda > 0$ and Condition 5.5 be satisfied. Define $u(x,y) = \lambda\overline{f}(x) - \underline{f}(y)$ and $v(x,y) = \lambda h(x) - h(y)$, where the $\overline{f}$ and $\underline{f}$ are the ones given in Theorem 5.2.*

*Then $u$ is a viscosity subsolution (in the sense of Definition 1.16) of*

$$(I - \alpha H_{2,\lambda})u = v.$$

PROOF. Let $\Phi_0 \in \mathcal{D}$ and $(\hat{x}, \hat{y}) \in E \times E$ satisfy

$$(u - \Phi_0)(\hat{x}, \hat{y}) = \sup_{x,y \in E} (u - \Phi_0)(x, y).$$

We want to show that

$$\alpha^{-1}(u(\hat{x}, \hat{y}) - v(\hat{x}, \hat{y})) \le (H_{2,\lambda}\Phi_0)^*(\hat{x}, \hat{y}). \tag{5.12}$$

We may assume $\Phi_0$ takes the following form:

$$\Phi_0(x, y) = \phi(x, y) + \rho(d_C(x, x_0) + d_C(y, y_0))$$

for some $x_0, y_0 \in E$, where

$$\phi(x, y) = \frac{\mu}{2}(\|x - y\|^2) + \frac{\gamma}{2}(\varphi(\|x\|^2) + \varphi(\|y\|^2))$$

takes the form in (5.9). Fix $\kappa > 0$; we let

$$\Phi(x, y) = \Phi_0(x, y) + \kappa(d_C(x, \hat{x}) + d_C(y, \hat{y})).$$

By Lemma 5.8,

$$(u - \Phi)(\hat{x}, \hat{y}) > (u - \Phi)(x, y) \qquad \forall\, (x, y) \neq (\hat{x}, \hat{y}).$$

We define

$$\Psi(x, y) = u(x, y) - \Phi(x, y),$$

$$\Psi_\varepsilon(x, y, \xi, \eta) = u(x, y) - \Phi(\xi, \eta) - \frac{1}{2\varepsilon}(\|x - \xi\|^2 + \|y - \eta\|^2)$$

and

$$\Psi_{\varepsilon,\delta}(x, y, \xi, \eta) = u(x, y) - \Phi(\xi, \eta) - \frac{1}{2\varepsilon}(\|x - \xi\|^2 + \|y - \eta\|^2) - \delta(\|C^0\xi\| + \|C^0\eta\|),$$

where $\|C^0\xi\| = +\infty$, if $\xi \notin \mathcal{D}(C)$.

We write the maximum of each of these functions:

$$M = \sup_{x,y \in E} \Psi(x, y) = \Psi(\hat{x}, \hat{y}),$$

$$M_\varepsilon = \sup_{x,y,\xi,\eta \in E} \Psi_\varepsilon(x, y, \xi, \eta),$$

$$M_{\varepsilon,\delta} = \sup_{x,y,\xi,\eta \in E} \Psi_{\varepsilon,\delta}(x, y, \xi, \eta).$$

It follows that $M \le M_\varepsilon$, $M_\varepsilon \ge M_{\varepsilon,\delta}$ and $M_\varepsilon \downarrow M$, as $\varepsilon \downarrow 0$ and $M_{\varepsilon,\delta} \uparrow M_\varepsilon$ as $\delta \downarrow 0$. See, for example, page 83 of [8]. The definitions of $M_\varepsilon, M_{\varepsilon,\delta}$ in [8] are slightly different than here, in the sense that suprema are taken locally for a ball of size $2r$ with arbitrary $r > 0$, instead of over the whole space.

Note that $\sup_{x,y}|u(x,y)|<+\infty$; note also that the form of $\phi$ in (5.9) implies that $\lim_{\|x\|+\|y\|\to+\infty}\Phi(x,y)=+\infty$, hence the suprema over the whole space are equal to the corresponding suprema over a sufficiently large open ball. Therefore, the same proof in [8] still works here.

(1) For each $\varepsilon,\theta>0$, we can choose $\delta=\delta(\varepsilon,\theta)$ and $(x_{\varepsilon,\theta},y_{\varepsilon,\theta},\xi_{\varepsilon,\theta},\eta_{\varepsilon,\theta})$ such that $\delta(\varepsilon,\theta)\downarrow 0$ as $\varepsilon,\theta\downarrow 0$ and

$$M_\varepsilon-\theta\le\Psi_{\varepsilon,\delta(\varepsilon,\theta)}(x_{\varepsilon,\theta},y_{\varepsilon,\theta},\xi_{\varepsilon,\theta},\eta_{\varepsilon,\theta}). \tag{5.13}$$

By Lemma 5.3, we can always select the $(x_{\varepsilon,\theta},y_{\varepsilon,\theta},\xi_{\varepsilon,\theta},\eta_{\varepsilon,\theta})$ so that (note that $\|C^0x\|$ as a function in $x$ is lower semicontinuous on $E$; see, e.g., Lemma 2.18 in [26])

$$\begin{aligned}&\Psi_{\varepsilon,\delta(\varepsilon,\theta)}(x,y,\xi,\eta)\\&\quad-\varepsilon(d_C(x,x_{\varepsilon,\theta})+d_C(y,y_{\varepsilon,\theta})+d_C(\xi,\xi_{\varepsilon,\theta})+d_C(\eta,\eta_{\varepsilon,\theta}))\\&\qquad\le\Psi_{\varepsilon,\delta(\varepsilon,\theta)}(x_{\varepsilon,\theta},y_{\varepsilon,\theta},\xi_{\varepsilon,\theta},\eta_{\varepsilon,\theta})\qquad\forall x,y,\xi,\eta\in E.\end{aligned} \tag{5.14}$$

Let $x=x_{\varepsilon,\theta},y=y_{\varepsilon,\theta}$; then

$$\begin{aligned}&\Phi(\xi_{\varepsilon,\theta},\eta_{\varepsilon,\theta})+\frac{1}{2\varepsilon}(\|x_{\varepsilon,\theta}-\xi_{\varepsilon,\theta}\|^2+\|y_{\varepsilon,\theta}-\eta_{\varepsilon,\theta}\|^2)\\&\quad+\delta(\varepsilon,\theta)(\|C^0\xi_{\varepsilon,\theta}\|+\|C^0\eta_{\varepsilon,\theta}\|)\\&\qquad\le\Phi(\xi,\eta)+\frac{1}{2\varepsilon}(\|x_{\varepsilon,\theta}-\xi\|^2+\|y_{\varepsilon,\theta}-\eta\|^2)\\&\qquad\quad+\delta(\varepsilon,\theta)(\|C^0\xi\|+\|C^0\eta\|)+\varepsilon(d_C(\xi,\xi_{\varepsilon,\theta})+d_C(\eta,\eta_{\varepsilon,\theta})).\end{aligned} \tag{5.15}$$

From $\Psi_\varepsilon\ge\Psi_{\varepsilon,\delta(\varepsilon,\theta)}$ and (5.13),

$$\Psi_\varepsilon(x,y,\xi,\eta)\le\Psi_\varepsilon(x_{\varepsilon,\theta},y_{\varepsilon,\theta},\xi_{\varepsilon,\theta},\eta_{\varepsilon,\theta})+\theta.$$

Take $x=x_{\varepsilon,\theta},y=y_{\varepsilon,\theta}$, therefore

$$\begin{aligned}&\theta+\Phi(\xi,\eta)+\frac{1}{2\varepsilon}(\|x_{\varepsilon,\theta}-\xi\|^2+\|y_{\varepsilon,\theta}-\eta\|^2)\\&\qquad\ge\Phi(\xi_{\varepsilon,\theta},\eta_{\varepsilon,\theta})\\&\qquad\quad+\frac{1}{2\varepsilon}(\|x_{\varepsilon,\theta}-\xi_{\varepsilon,\theta}\|^2+\|y_{\varepsilon,\theta}-\eta_{\varepsilon,\theta}\|^2),\qquad\xi,\eta\in E.\end{aligned} \tag{5.16}$$

(2) It follows from (5.13) that $\lim_{\varepsilon\to0+,\theta\to0+}\Psi_{\varepsilon,\delta(\varepsilon,\theta)}(x_{\varepsilon,\theta},y_{\varepsilon,\theta})=M$, hence

$$\begin{aligned}&(x_{\varepsilon,\theta},y_{\varepsilon,\theta}),(\xi_{\varepsilon,\theta},\eta_{\varepsilon,\theta})\to(\hat{x},\hat{y}),\\&\qquad\overline{f}(x_{\varepsilon,\theta})\to\overline{f}(\hat{x}),\qquad\underline{f}(x_{\varepsilon,\theta})\to\underline{f}(\hat{y})\qquad\text{as }\varepsilon,\theta\downarrow0\end{aligned} \tag{5.17}$$

(e.g., Step 3 on page 84 of [8] and Lemma 5.8).

Take $\xi = \xi_{\varepsilon,\theta}$, $\eta = \eta_{\varepsilon,\theta}$, $y = y_{\varepsilon,\theta}$ in (5.14):

$$x \to \lambda \overline{f}(x) - \frac{1}{2\varepsilon}\|x - \xi_{\varepsilon,\theta}\|^2 - \varepsilon d_C(x, x_{\varepsilon,\theta})$$

has a maximum at $x_{\varepsilon,\theta}$. Since $\overline{f}$ is a subsolution of (5.2) in the sense of Definition 1.16,

$$\begin{aligned}\lambda\frac{\overline{f} - h}{\alpha}(x_{\varepsilon,\theta}) \le \frac{1}{\varepsilon}\langle x_{\varepsilon,\theta} - \xi_{\varepsilon,\theta}, C^0\xi_{\varepsilon,\theta}\rangle + \varepsilon \\ + \sup_{\|p\|\le\varepsilon} \lambda G\Big(x_{\varepsilon,\theta}, \frac{1}{\lambda}\Big(\frac{x_{\varepsilon,\theta} - \xi_{\varepsilon,\theta}}{\varepsilon} + p\Big)\Big).\end{aligned}$$

Similarly, noting $y_{\varepsilon,\theta}$ is a maximum of

$$y \to -\underline{f}(y) - \frac{1}{2\varepsilon}\|y - \eta_{\varepsilon,\theta}\|^2 - \varepsilon d_C(y, y_{\varepsilon,\theta}),$$

by the supersolution property of $\underline{f}$,

$$\begin{aligned}\frac{\underline{f} - h}{\alpha}(y_{\varepsilon,\theta}) \ge -\frac{1}{\varepsilon}\langle y_{\varepsilon,\theta} - \eta_{\varepsilon,\theta}, C^0\eta_{\varepsilon,\theta}\rangle - \varepsilon \\ + \inf_{\|q\|\le\varepsilon} G\Big(y_{\varepsilon,\theta}, -\frac{y_{\varepsilon,\theta} - \eta_{\varepsilon,\theta}}{\varepsilon} + q\Big).\end{aligned}$$

Therefore

$$\begin{aligned}&\alpha^{-1}(u - v)(x_{\varepsilon,\theta}, y_{\varepsilon,\theta}) \\ &\quad\le \Big\langle \frac{1}{\varepsilon}(x_{\varepsilon,\theta} - \xi_{\varepsilon,\theta}, y_{\varepsilon,\theta} - \eta_{\varepsilon,\theta}), (C^0\xi_{\varepsilon,\theta}, C^0\eta_{\varepsilon,\theta})\Big\rangle + 2\varepsilon \\ &\quad\quad + \sup_{\|p\|^2+\|q\|^2\le 2\varepsilon^2}\Big(\lambda G\Big(x_{\varepsilon,\theta}, \frac{1}{\lambda}\Big(\frac{x_{\varepsilon,\theta} - \xi_{\varepsilon,\theta}}{\varepsilon} + p\Big)\Big) \\ &\quad\quad\quad - G\Big(y_{\varepsilon,\theta}, -\frac{y_{\varepsilon,\theta} - \eta_{\varepsilon,\theta}}{\varepsilon} + q\Big)\Big).\end{aligned} \tag{5.18}$$

(3) Apply Lemma A.8 of [8] to (5.16), and take $\theta = \varepsilon^3$. We denote

$$x_\varepsilon = x_{\varepsilon,\varepsilon^3}, \qquad y_\varepsilon = y_{\varepsilon,\varepsilon^3}, \qquad \xi_\varepsilon = \xi_{\varepsilon,\varepsilon^3}, \qquad \eta_\varepsilon = \eta_{\varepsilon,\varepsilon^3}.$$

(a) After some algebra (for details, see page 86 of [8]),

$$\limsup_{\varepsilon\to 0+}\Big\| D\phi(\hat{x}, \hat{y}) - \Big(\frac{x_\varepsilon - \xi_\varepsilon}{\varepsilon}, \frac{y_\varepsilon - \eta_\varepsilon}{\varepsilon}\Big)\Big\| \le \rho + \kappa. \tag{5.19}$$

(b) Take $\theta = \varepsilon^3, \xi = S(h)\xi_\varepsilon, \eta = S(h)\eta_\varepsilon$ in (5.15); noting $\|C^0 S(h)\xi_\varepsilon\| \leq \|C^0\xi_\varepsilon\|$ and $\|C^0 S(h)\eta_\varepsilon\| \leq \|C^0\eta_\varepsilon\|$,

$$\begin{aligned}\frac{1}{2\varepsilon}\frac{1}{h}(\|(x_\varepsilon, y_\varepsilon) - (\xi_\varepsilon, \eta_\varepsilon)\|^2 - \|(x_\varepsilon, y_\varepsilon) - \mathcal{S}(h)(\xi_\varepsilon, \eta_\varepsilon)\|^2)&\\ \leq \frac{1}{h}(\Phi(S(h)\xi_\varepsilon, S(h)\eta_\varepsilon) - \Phi(\xi_\varepsilon, \eta_\varepsilon))&\\ + \varepsilon\frac{1}{h}(d_C(S(h)\xi_\varepsilon, \xi_\varepsilon) + d_C(S(h)\eta_\varepsilon, \eta_\varepsilon)).&\end{aligned}$$

Send $h \to 0+$; by (3.2) and (3.4),

$$\begin{aligned}\left\langle \frac{(x_\varepsilon, y_\varepsilon) - (\xi_\varepsilon, \eta_\varepsilon)}{\varepsilon}, (C^0\xi_\varepsilon, C^0\eta_\varepsilon)\right\rangle&\\ \leq \limsup_{h\to 0+} \frac{\Phi(S(h)\xi_\varepsilon, S(h)\eta_\varepsilon) - \Phi(\xi_\varepsilon, \eta_\varepsilon)}{h} + 2\varepsilon.&\end{aligned}$$

Hence by (5.17),

$$\begin{aligned}\lim_{\varepsilon\to 0+}\left\langle \frac{(x_\varepsilon, y_\varepsilon) - (\xi_\varepsilon, \eta_\varepsilon)}{\varepsilon}, (C^0\xi_\varepsilon, C^0\eta_\varepsilon)\right\rangle&\\ \leq \limsup_{h\to 0+,\varepsilon\to 0+} \frac{\Phi(S(h)\xi_\varepsilon, S(h)\eta_\varepsilon) - \Phi(\xi_\varepsilon, \eta_\varepsilon)}{h}&\\ \leq \limsup_{h\to 0+,(x,y)\to(\hat{x},\hat{y})} \frac{\Phi(S(h)x, S(h)y) - \Phi(x,y)}{h}&\\ \leq D^+_{\mathcal{C}}\Phi(\hat{x},\hat{y}) \leq D^+_{\mathcal{C}}\Phi_0(\hat{x},\hat{y}) + 2\kappa.&\end{aligned} \tag{5.20}$$

Finally, apply (5.19) and (5.20) to (5.18); noting (5.17),

$$\begin{aligned}\frac{u-v}{\alpha}(\hat{x},\hat{y}) \leq D^+_{\mathcal{C}}\Phi_0(\hat{x},\hat{y}) + 2\kappa&\\ + \sup_{\|p\|^2+\|q\|^2\leq 2(\rho+\kappa)^2} G_\lambda((\hat{x},\hat{y}); D\phi(\hat{x},\hat{y}) + (p,q)).&\end{aligned} \tag{5.21}$$

Send $\kappa \to 0$;

$$\frac{u-v}{\alpha}(\hat{x},\hat{y}) \leq (H_{2,\lambda}\Phi_0)^*(\hat{x},\hat{y}).$$ □

We discuss a property of $D^+_{\mathcal{C}}$.

Lemma 5.10. *Let $\Phi(x,y) \in \mathcal{D}$ be such that*

$$\Phi(x,y) = \phi(x,y) + \theta(d_C(x,x_0) + d_C(y,y_0)), \qquad x_0, y_0 \in E, \theta > 0,$$

*where*

$$\phi(x,y)=\frac{\mu}{2}\|x-y\|^2+\frac{\varepsilon}{2}(\varphi(\|x\|^2)+\varphi(\|y\|^2)),$$
$$\mu>0,\varphi\in C^1([0,\infty)),\varphi'\geq 0.$$

*Then*

$$D_{\mathcal{C}}^{+}\phi(x,y)\leq 0 \tag{5.22}$$

*and*

$$D_{\mathcal{C}}^{+}\Phi(x,y)\leq 2\theta. \tag{5.23}$$

PROOF. Equation (5.22) follows because $S(t)$ is a contraction semigroup with the property $\|S(t)x\|\leq\|x\|$:

$$\phi(S(h)x,S(h)y)\leq\phi(x,y)\qquad\forall\, x,y\in E,h>0.$$

For each $x,x_0\in E$, by the definition of $d_C$, there exists $t_0>0$ such that

$$\begin{aligned}d_C(x,x_0)+h&=\inf_{t\geq 0}(t+\|x-S(t)x_0\|)+h=t_0+\|x-S(t_0)x_0\|+h\\&\geq t_0+h+\|S(h)x-S(t_0+h)x_0\|\geq d_C(S(h)x,x_0)\qquad\forall\, h>0.\end{aligned}$$

Hence (5.23) follows. □

PROOF OF THEOREM 5.2. We assume $(\overline{f}-\underline{f})(z_0)=4\delta_0>0$ for some $z_0\in E$ (otherwise, there is nothing to prove). We want to create a contradiction.

Let $1<\lambda<1+\delta_0/(1+\sup_x|\overline{f}(x)|)$. We recall that $u(x,y)=\lambda\overline{f}(x)-\underline{f}(y)$. Then $u(z_0,z_0)\geq 3\delta_0>0$. Let $\varphi(r)$ be given by Condition 5.5(2); we define

$$\phi(x,y)=\frac{\mu}{2}\|x-y\|^2+\frac{\varepsilon}{2}(\varphi(\|x\|^2)+\varphi(\|y\|^2)).$$

Hence for $0<\varepsilon<(\delta_0)/(1+\varphi(|z_0|^2))$,

$$0<2\delta_0\leq u(z_0,z_0)-\delta_0\leq(u-\phi)(z_0,z_0)\leq\sup_{x,y\in E}(u-\phi)(x,y). \tag{5.24}$$

By Lemma 5.3 on perturbed optimization, for each $\theta>0$, there exist $x_{\mu,\varepsilon,\theta},y_{\mu,\varepsilon,\theta}\in E$ such that

$$(u-\phi)(x,y)-\theta(d_C(x,x_{\mu,\varepsilon,\theta})+d_C(y,y_{\mu,\varepsilon,\theta}))\leq(u-\phi)(x_{\mu,\varepsilon,\theta},y_{\mu,\varepsilon,\theta})$$

and

$$\sup_{x,y\in E}(u-\phi)(x,y)\leq(u-\phi)(x_{\mu,\varepsilon,\theta},y_{\mu,\varepsilon,\theta})+\theta. \tag{5.25}$$

Denote

$$\Phi(x,y)=\phi(x,y)+\theta(d_C(x,x_{\mu,\varepsilon,\theta})+d_C(y,y_{\mu,\varepsilon,\theta})),$$

therefore $u(x,y)-\Phi(x,y)$ attains its maximum at $(x_{\mu,\varepsilon,\theta}, y_{\mu,\varepsilon,\theta})$. From (5.24) and (5.25),

$$\begin{aligned} \delta_0 &\le \delta_0 + \frac{\mu}{2}\|x_{\mu,\varepsilon,\theta} - y_{\mu,\varepsilon,\theta}\|^2 + \frac{\varepsilon}{2}(\varphi(\|x_{\mu,\varepsilon,\theta}\|^2) + \varphi(\|y_{\mu,\varepsilon,\theta}\|^2)) \\ &\le u(x_{\mu,\varepsilon,\theta}, y_{\mu,\varepsilon,\theta}) \le \lambda \sup_x |\overline{f}(x)| + \sup_x |\underline{f}(x)| < \infty \end{aligned} \tag{5.26}$$

for every $\theta \le \delta_0, |\lambda - 1| \le \delta_0(1 + \sup_x |\overline{f}(x)|)^{-1}$. Equation (5.26) implies the existence of constants $C_\mu, M_\varepsilon$ such that

$$\mu\|x_{\mu,\varepsilon,\theta} - y_{\mu,\varepsilon,\theta}\| \le C_\mu < \infty \qquad \forall\, \varepsilon > 0, 0 < \theta < \delta_0,$$

$$\sqrt{\|x_{\mu,\varepsilon,\theta}\|^2 + \|y_{\mu,\varepsilon,\theta}\|^2} \le M_\varepsilon < \infty \qquad \forall\, \mu > 0, 0 < \theta < \delta_0.$$

In addition, since

$$C_0 \equiv \sup_{r\ge 0} |r\varphi'(r^2)| + |r\varphi'(r)| < \infty,$$

there exists constant $N_{\varepsilon,\mu}$,

$$\begin{aligned} &\|x_{\mu,\varepsilon,\theta}\| + \|y_{\mu,\varepsilon,\theta}\| + \mu\|x_{\mu,\varepsilon,\theta} - y_{\mu,\varepsilon,\theta}\| \\ &\qquad + \varepsilon\varphi'(\|x_{\mu,\varepsilon,\theta}\|^2)\|x_{\mu,\varepsilon,\theta}\| + \varepsilon\varphi'(\|y_{\mu,\varepsilon,\theta}\|^2)\|y_{\mu,\varepsilon,\theta}\| \le N_{\varepsilon,\mu} < \infty, \end{aligned}$$

for $0 < \theta < \delta_0$. Since

$$\begin{aligned} &D\phi(x_{\mu,\varepsilon,\theta}, y_{\mu,\varepsilon,\theta}) \\ &\qquad = (\mu(x_{\mu,\varepsilon,\theta} - y_{\mu,\varepsilon,\theta}) + \varepsilon\varphi'(\|x_{\mu,\varepsilon,\theta}\|^2)x_{\mu,\varepsilon,\theta}, \\ &\qquad\quad - \mu(x_{\mu,\varepsilon,\theta} - y_{\mu,\varepsilon,\theta}) + \varepsilon\varphi'(\|y_{\mu,\varepsilon,\theta}\|^2)y_{\mu,\varepsilon,\theta}), \end{aligned}$$

and $u$ is a viscosity subsolution of $(I - \alpha H_{2,\lambda})u = v$ (Lemma 5.9), by the estimate in (5.23) and (5.26), for $\theta < \delta_0$ and $|\lambda - 1| \le \delta_0(1 + 2\sup_x |h(x)| \vee 2\sup_x |\overline{f}(x)|)^{-1}$,

$$\begin{aligned} &\frac{\delta_0}{2} - (h(x_{\mu,\varepsilon,\theta}) - h(y_{\mu,\varepsilon,\theta})) \\ &\qquad \le u(x_{\mu,\varepsilon,\theta}, y_{\mu,\varepsilon,\theta}) - (\lambda h(x_{\mu,\varepsilon,\theta}) - h(y_{\mu,\varepsilon,\theta})) \\ &\qquad \le 2\theta + \sup_{\|p\|^2 + \|q\|^2 \le 2\theta^2} G_\lambda((x_{\mu,\varepsilon,\theta}, y_{\mu,\varepsilon,\theta}); \\ &\qquad\qquad\qquad\qquad\qquad D\phi(x_{\mu,\varepsilon,\theta}, y_{\mu,\varepsilon,\theta}) + (p,q)) \\ &\qquad \le 2\theta + \sup_{\|p\|^2 + \|q\|^2 \le 2\theta^2} G_\lambda\{(x_{\mu,\varepsilon,\theta}, y_{\mu,\varepsilon,\theta}); \\ &\qquad\qquad\qquad\qquad\qquad (\mu(x_{\mu,\varepsilon,\theta} - y_{\mu,\varepsilon,\theta}) \\ &\qquad\qquad\qquad\qquad\qquad\ + \varepsilon\varphi'(\|x_{\mu,\varepsilon,\theta}\|^2)x_{\mu,\varepsilon,\theta} + p, \end{aligned} \tag{5.27}$$

$$- \mu(x_{\mu,\varepsilon,\theta} - y_{\mu,\varepsilon,\theta}) + \varepsilon\varphi'(\|y_{\mu,\varepsilon,\theta}\|^2)y_{\mu,\varepsilon,\theta} + q)\}.$$

We select $\lambda_0, \lambda_1$ satisfying $1 < \lambda_0 < \lambda_1 < \lambda$, and let them be fixed. By Condition 5.5, for $\|p\| < 1$,

$$\begin{aligned}
&\lambda G\Big(x_{\mu,\varepsilon,\theta}, \frac{\mu(x_{\mu,\varepsilon,\theta} - y_{\mu,\varepsilon,\theta}) + \varepsilon\varphi'(\|x_{\mu,\varepsilon,\theta}\|^2)x_{\mu,\varepsilon,\theta} + p}{\lambda}\Big) \\
&\quad = \lambda G\Big(x_{\mu,\varepsilon,\theta}, \frac{\mu(x_{\mu,\varepsilon,\theta} - y_{\mu,\varepsilon,\theta}) + \varepsilon\varphi'(\|x_{\mu,\varepsilon,\theta}\|^2)x_{\mu,\varepsilon,\theta} + p}{\lambda}\Big) \\
&\qquad - \lambda G\Big(x_{\mu,\varepsilon,\theta}, \frac{\mu(x_{\mu,\varepsilon,\theta} - y_{\mu,\varepsilon,\theta}) + \varepsilon\varphi'(\|x_{\mu,\varepsilon,\theta}\|^2)x_{\mu,\varepsilon,\theta}}{\lambda}\Big) \\
&\qquad + \lambda G\Big(x_{\mu,\varepsilon,\theta}, \frac{\mu(x_{\mu,\varepsilon,\theta} - y_{\mu,\varepsilon,\theta}) + \varepsilon\varphi'(\|x_{\mu,\varepsilon,\theta}\|^2)x_{\mu,\varepsilon,\theta}}{\lambda}\Big) \\
&\qquad - \lambda_1 G\Big(x_{\mu,\varepsilon,\theta}, \frac{\mu(x_{\mu,\varepsilon,\theta} - y_{\mu,\varepsilon,\theta})}{\lambda_1}\Big) + \lambda_1 G\Big(x_{\mu,\varepsilon,\theta}, \frac{\mu(x_{\mu,\varepsilon,\theta} - y_{\mu,\varepsilon,\theta})}{\lambda_1}\Big) \\
&\quad \le \lambda\sigma_{N_{\varepsilon,\mu}+1}\Big(\frac{\|p\|}{\lambda}\Big) + \lambda_1\gamma_{\lambda/\lambda_1, C_\mu/\lambda_1}\Big(\frac{\varepsilon}{\lambda_1}\Big) \\
&\qquad + \lambda_1 G\Big(x_{\mu,\varepsilon,\theta}, \frac{\mu(x_{\mu,\varepsilon,\theta} - y_{\mu,\varepsilon,\theta})}{\lambda_1}\Big).
\end{aligned}$$

Similarly, for $\|q\| < 1$ and $0 \le \varepsilon < 1$,

$$\begin{aligned}
&G(y_{\mu,\varepsilon,\theta}, \mu(x_{\mu,\varepsilon,\theta} - y_{\mu,\varepsilon,\theta}) - \varepsilon\varphi'(\|y_{\mu,\varepsilon,\theta}\|^2)y_{\mu,\varepsilon,\theta} - q) \\
&\quad \ge G(y_{\mu,\varepsilon,\theta}, \mu(x_{\mu,\varepsilon,\theta} - y_{\mu,\varepsilon,\theta}) - \varepsilon\varphi'(\|y_{\mu,\varepsilon,\theta}\|^2)y_{\mu,\varepsilon,\theta}) - \sigma_{N_{\varepsilon,\mu}+1}(\|q\|) \\
&\quad \ge \lambda_0 G\Big(y_{\mu,\varepsilon,\theta}, \frac{\mu(x_{\mu,\varepsilon,\theta} - y_{\mu,\varepsilon,\theta})}{\lambda_0}\Big) - \gamma_{\lambda_0, C_\mu + C_0}(\varepsilon) - \sigma_{N_{\varepsilon,\mu}+1}(\|q\|).
\end{aligned}$$

Therefore

$$\begin{aligned}
&G_\lambda\{(x_{\mu,\varepsilon,\theta}, y_{\mu,\varepsilon,\theta}); (\mu(x_{\mu,\varepsilon,\theta} - y_{\mu,\varepsilon,\theta}) + \varepsilon\varphi'(\|x_{\mu,\varepsilon,\theta}\|^2)x_{\mu,\varepsilon,\theta} + p, \\
&\qquad\qquad - \mu(x_{\mu,\varepsilon,\theta} - y_{\mu,\varepsilon,\theta}) + \varepsilon\varphi'(\|y_{\mu,\varepsilon,\theta}\|^2)y_{\mu,\varepsilon,\theta} + q)\} \\
&\quad \le \lambda_1 G\Big(x_{\mu,\varepsilon,\theta}, \frac{\mu(x_{\mu,\varepsilon,\theta} - y_{\mu,\varepsilon,\theta})}{\lambda_1}\Big) - \lambda_0 G\Big(y_{\mu,\varepsilon,\theta}, \frac{\mu(x_{\mu,\varepsilon,\theta} - y_{\mu,\varepsilon,\theta})}{\lambda_0}\Big) \\
&\qquad + \lambda\sigma_{N_{\varepsilon,\mu}+1}\Big(\frac{\|p\|}{\lambda}\Big) + \lambda_1\gamma_{\lambda/\lambda_1, C_\mu/\lambda_1}\Big(\frac{\varepsilon}{\lambda_1}\Big) \\
&\qquad + \gamma_{\lambda_0, C_\mu + C_0}(\varepsilon) + \sigma_{N_{\varepsilon,\mu}+1}(\|q\|) \\
&\quad \le \lambda_0\rho_{\lambda_1/\lambda_0}\Big(\|x_{\mu,\varepsilon,\theta} - y_{\mu,\varepsilon,\theta}\| + \frac{\mu}{\lambda_0}\|x_{\mu,\varepsilon,\theta} - y_{\mu,\varepsilon,\theta}\|^2\Big) + \lambda\sigma_{N_{\varepsilon,\mu}+1}\Big(\frac{\|p\|}{\lambda}\Big)
\end{aligned}$$

$$+\lambda_1\gamma_{\lambda/\lambda_1,C_\mu/\lambda_1}\Big(\frac{\varepsilon}{\lambda_1}\Big)+\gamma_{\lambda_0,C_\mu+C_0}(\varepsilon)+\sigma_{N_\varepsilon,\mu+1}(\|q\|).$$

We rewrite (5.27) next. For $\theta<\delta_0<1$ and $1<\lambda\le 1+\frac{\delta_0}{1+2\|h\|\vee 2\|\overline{f}\|}$,

$$\begin{aligned}&\frac{\delta_0}{2}-(h(x_{\mu,\varepsilon,\theta})-h(y_{\mu,\varepsilon,\theta}))\\ &\qquad\le 2\theta+\sup_{\|p\|^2+\|q\|^2\le 2\theta^2}\Big(\lambda\sigma_{N_\varepsilon,\mu+1}\Big(\frac{\|p\|}{\lambda}\Big)+\sigma_{N_\varepsilon,\mu+1}(\|q\|)\Big)\\ &\qquad\quad+\lambda_1\gamma_{\lambda/\lambda_1,C_\mu/\lambda_1}\Big(\frac{\varepsilon}{\lambda_1}\Big)+\gamma_{\lambda_0,C_\mu+C_0}(\varepsilon)\\ &\qquad\quad+\lambda_0\rho_{\lambda_1/\lambda_0}\Big(\|x_{\mu,\varepsilon,\theta}-y_{\mu,\varepsilon,\theta}\|+\frac{\mu}{\lambda_0}\|x_{\mu,\varepsilon,\theta}-y_{\mu,\varepsilon,\theta}\|^2\Big).\end{aligned}\tag{5.28}$$

Let

$$m_\mu=\sup_{x,y\in E}\Big(u(x,y)-\frac{\mu}{2}\|x-y\|^2\Big).$$

From (5.25),

$$\begin{aligned}u(x,y)-\frac{\mu}{2}\|x-y\|^2&\le\liminf_{\varepsilon\to0,\theta\to0}\Big(u(x_{\mu,\varepsilon,\theta},y_{\mu,\varepsilon,\theta})-\frac{\mu}{2}\|x_{\mu,\varepsilon,\theta}-y_{\mu,\varepsilon,\theta}\|^2\Big)\\ &\le\limsup_{\varepsilon\to0,\theta\to0}\Big(u(x_{\mu,\varepsilon,\theta},y_{\mu,\varepsilon,\theta})-\frac{\mu}{2}\|x_{\mu,\varepsilon,\theta}-y_{\mu,\varepsilon,\theta}\|^2\Big)\\ &\le m_\mu,\qquad x,y\in E.\end{aligned}$$

Hence

$$m_\mu=\lim_{\varepsilon\to0,\theta\to0}\Big(u(x_{\mu,\varepsilon,\theta},y_{\mu,\varepsilon,\theta})-\frac{\mu}{2}\|x_{\mu,\varepsilon,\theta}-y_{\mu,\varepsilon,\theta}\|^2\Big).$$

Appling Lemma 3.2 in [8],

$$\lim_{\mu\to\infty}\limsup_{\varepsilon\to0,\theta\to0}\mu\|x_{\mu,\varepsilon,\theta}-y_{\mu,\varepsilon,\theta}\|^2=0.$$

In (5.28), let $\theta\downarrow0$, then $\varepsilon\downarrow0$, then $\mu\uparrow+\infty$; we obtain

$$0<\delta_0/2\le 0.$$

A contradiction. □

## APPENDIX

### A.1. Verifying semigroup generation condition.

A.1.1. *The Cahn–Hilliard equation.* Let $\omega \equiv \sup_r |V''(r)|^2/4$ and let $C, C_n$ be defined according to (4.4) and (1.23). We prove the following.

LEMMA A.1. *The closure of $C$ (resp. $C_n$) is an $m$-dissipative operator in $L^2(\mathcal{O})$.*

The proof is divided into two parts.

LEMMA A.2. *Both $C$ and $C_n$ are dissipative.*

PROOF. Let $x, y \in \mathcal{D}(C)$; then

$$\begin{aligned}
&\langle Cx - Cy, x - y\rangle \\
&\qquad = \langle \Delta(-\Delta(x-y)), x-y\rangle + \langle V'(x) - V'(y), \Delta(x-y)\rangle - \omega\|x-y\|^2 \\
&\qquad \le -\|\Delta(x-y)\|^2 + \sup_r |V''(r)|\|x-y\|\|\Delta(x-y)\| \\
&\qquad\quad - \frac{\sup_r |V''(r)|^2}{4}\|x-y\|^2 \le 0.
\end{aligned}$$

The case of $C_n$ can be treated similarly. □

LEMMA A.3. *Let $0 < \alpha \le \alpha_0$ where $\alpha_0$ is some prefixed small number. Suppose $y_0 \in L^2(\mathcal{O})$. Then for each $y_n \equiv P_n y_0$, there exists $x_n \in \mathcal{R}(P_n)$ such that*

$$(I - \alpha C_n)x_n = y_n. \tag{A.1}$$

*Moreover,*

$$\|x_n\| \le \|y_n\| \quad \text{and} \quad \|C_n x_n\| \le (2/\alpha)\|y_n\|, \tag{A.2}$$

$$\|\Delta x_n\| \le c(1 + \|y_n\|), \tag{A.3}$$

*where $c$ is a constant depending on $F$ and $\alpha$. Consequently, there exists $x_0 \in \mathcal{D}(\bar{C})$ such that*

$$(I - \alpha\bar{C})x_0 = y_0.$$

PROOF. For each $n$ fixed and finite, by its definition (1.23), $C_n$ is Lipschitz on Range$(P_n)$. By the fixed point theorem, the essentially finite-dimensional equation (A.1) has a solution when $\alpha > 0$ is sufficiently small. Since $C_n$ is dissipative on Range$(P_n)$, solution actually exists for all $\alpha > 0$. See Lemma 2.13 of [26].

$\|x_n\| \le \|y_n\|$ follows from the dissipativity of $C_n$. It follows then that $\alpha\|C_n x_n\| = \|x_n - y_n\| \le 2\|y_n\|$. $\langle x_n - \alpha C_n x_n, x_n\rangle = \langle y_n, x_n\rangle$. That is,

$$\|x_n\|^2 - \alpha(-\|\Delta x_n\|^2 + \langle V'(x_n), \Delta x_n\rangle - \omega\|x_n\|^2) = \langle y_n, x_n\rangle,$$

which implies

$$\alpha(\|\Delta x_n\|^2 + \omega\|x_n\|^2) \le \|y_n\|\|x_n\| - \|x_n\|^2 + \alpha\|V'(x_n)\|\|\Delta x_n\|.$$

Since $V'(r)$ grows at most linearly, we can find constant $c > 0$,

$$\|V'(x_n)\| \le c(1 + \|x_n\|).$$

Hence

$$\begin{aligned}\alpha\|\Delta x_n\|^2 &\le \|y_n\|\|x_n\| + \alpha c(1 + \|x_n\|)\|\Delta x_n\| \\ &\le \|y_n\|^2 + \alpha c(1 + \|y_n\|)\|\Delta x_n\|,\end{aligned}$$

implying (A.3).

For each $x \in H^2(\mathcal{O})$, $\|\Delta x\| < \infty$ and

$$\Delta V'(x) = V''(x)\Delta x + V'''(x)\nabla x \nabla x. \tag{A.4}$$

We note that the Sobolev embedding $H^1(\mathcal{O}) \to L^4(\mathcal{O})$ holds for space dimensions $d = 1, 2, 3$. Such result can be found in [1]: the case of $d = 3$ follows from Lemma 5.10, the case of $d = 2$ from Corollary 5.13 and the case of $d = 1$ follows from Corollary 5.16 of [1]. Therefore,

$$\begin{aligned}\|\Delta V'(x_n)\| &\le \sup_r |V''(r)| \|\Delta x_n\|_{L^2(\mathcal{O})} + \sup_r |V'''(r)| \|\nabla x_n\|^2_{L^4(\mathcal{O})} \\ &\le C(1 + \|\Delta x_n\|^2_{L^2(\mathcal{O})})\end{aligned}$$

for some constant $C$ independent of the $x_n$'s. By (A.3), $\sup_n \|\Delta V'(x_n)\| < \infty$. Using this estimate and $\sup_n \|C_n x_n\| < \infty$, we obtain $\sup_n \|\Delta^2 x_n\| < \infty$.

The boundedness of $\sup_n(\|\Delta x_n\| + \|\Delta^2 x_n\|)$ implies that $\Delta x_n$ is relatively compact in $L^2(\mathcal{O})$. Similarly, the boundedness of $\|\Delta x_n\|$ and $\|x_n\|$ implies the relative compactness of $\nabla x_n$ and $x_n$. Selecting a subsequence if necessary, we have $x_n \to x_0, \Delta x_n \to \Delta x_0$ and $\nabla x_n \to \nabla x_0$ for some $x_0 \in H^2(\mathcal{O})$. By (A.4), $\Delta V'(x_n) \to \Delta V'(x_0)$. Therefore

$$\|C_n x_n - C x_n\| = \|(P_n - I)\Delta V'(x_n)\| \to 0.$$

Noting

$$\alpha C x_n = \alpha C_n x_n + \alpha(C x_n - C_n x_n) = x_n - y_n + \alpha(C x_n - C_n x_n) \to x_0 - y_0,$$

$(x_0, \alpha^{-1}(x_0 - y_0)) \in \bar{C}$. By (A.1),

$$y_0 = x_0 - \alpha \frac{x_0 - y_0}{\alpha} \in (I - \alpha\bar{C})x_0. \qquad \square$$

A.1.2. *The quasilinear equation with viscosity.* We consider the $C_n, C$ in (1.27) and (4.5). Using similar a priori estimate arguments as in the Cahn–Hilliard equation case, we can prove the following:

LEMMA A.4. *$C_n$ and the closure of $C$ are both $m$-dissipative operators in $L^2(\mathcal{O})$.*

A.1.3. *The Allen–Cahn equation.* Let $C_n$ and $C$ be defined according to (1.18) and (4.1). Again, using a priori estimate arguments similar to the Cahn–Hilliard case, we can prove that both $C_n$ and $C$ are $m$-dissipative operators in $L^2(\mathcal{O})$. Alternatively, this conclusion can also be established by invoking the classical perturbation theory in, for instance, Corollary 6.19(i) of [26].

**A.2. Exponential compact containment estimates.** We illustrate the use of a stochastic Lyapunov function technique to verify Condition 1.12. We consider Examples 1.2, 1.5 and 1.8.

A.2.1. *Stochastic Allen–Cahn equation.* Recall (1.17) in Example 1.2:

$$dX_n(t) = \Delta P_n X_n(t)\,dt - P_n V'(P_n X_n(t))\,dt + \frac{1}{\sqrt{n}} B_n(X_n(t))\,dW(t). \tag{A.5}$$

The goal of this subsection is to prove the following lemma.

Lemma A.5. *Condition 1.12 holds for $X_n$.*

We introduce the *free energy* function

$$\mathcal{E}(x) \equiv \tfrac{1}{2}\|\nabla x\|^2 + \int_{\mathcal{O}} V(x)\,d\theta.$$

First, we prove the following estimate: for every $T, a > 0$ and $C_0 > 0$, there exists $C_1 > 0$ such that

$$\sup_{x : \mathcal{E}(x) \le C_0} P(\mathcal{E}(X_n(t)) > C_1, \text{ some } 0 < t \le T | X_n(0) = x) \le e^{-na}. \tag{A.6}$$

Let us approximate $\mathcal{E}$ by

$$\mathcal{E}_n(x) \equiv -\tfrac{1}{2}\langle \Delta P_n x, x\rangle + \int_{\mathcal{O}} V(P_n x(\theta))\,d\theta. \tag{A.7}$$

Note that if $X_n(0) \in \mathcal{R}(P_n)$, the range of $P_n$, then $X_n(t) \in \mathcal{R}(P_n)$, hence $\mathcal{E}_n(X_n(t)) = \mathcal{E}(X_n(t))$. Define

$$f_n(x) \equiv \log\left(1 + \frac{1}{M^2}\mathcal{E}_n(x)\right) \tag{A.8}$$

where $M \equiv \sup_{\theta,x} |\sigma(\theta, x)| < \infty$. Then

$$Df_n(x) = \frac{(-1/M^2)(\Delta P_n x - P_n V'(P_n x))}{1 + (1/M^2)\mathcal{E}_n(x)}$$

and

$$D^2 f_n(x) = \frac{(-1/M^2)(\Delta P_n x - P_n V'(P_n x)) \otimes (1/M^2)(\Delta P_n x - P_n V'(P_n x))}{(1+(1/M^2)\mathcal{E}_n(x))^2} + \frac{(1/M^2)(-\Delta P_n + P_n V''(P_n x))}{1+(1/M^2)\mathcal{E}_n(x)},$$

where $P_n V''(P_n x)$ means a linear operator on $L^2(\mathcal{O})$ [for each $x \in L^2(\mathcal{O})$ fixed]:

$$(P_n V''(P_n x))y \equiv \sum_{k\equiv(k_1,\dots,k_d)=(1,\dots,1)}^{(m_n,\dots,m_n)} \langle V''(P_n x)y, e_k\rangle e_k \qquad \forall\, y \in L^2(\mathcal{O}).$$

Then

$$\begin{aligned} H_n f_n(x) &\equiv \langle \Delta P_n x - P_n V'(P_n x), Df_n(x)\rangle \\ &\quad + \frac{1}{2}\|B_n^*(x) Df_n(x)\|_{U_0}^2 + \frac{1}{2n}\operatorname{Tr}(D^2 f_n(x) B_n(x) B_n^*(x)) \\ &= \frac{(-1/M^2)\|\Delta P_n x - P_n V'(P_n x)\|^2}{1+(1/M^2)\mathcal{E}_n(x)} \\ &\quad + \frac{1}{2}\left(1-\frac{1}{n}\right)\left\|\frac{(1/M^2)B_n^*(x)(\Delta P_n x - P_n V_n'(P_n x))}{1+(1/M^2)\mathcal{E}_n(x)}\right\|^2 \\ &\quad + \frac{1}{2n}\left(\sum_{k=(1,\dots,1)}^{(m_n,\dots,m_n)}\left\langle (-\Delta P_n + P_n V''(P_n x))\frac{B_n(x)}{M^2}e_k, \frac{B_n(x)}{M^2}e_k\right\rangle\right) \\ &\quad \times \left(1+\frac{1}{M^2}\mathcal{E}_n(x)\right)^{-1} \\ &= \frac{(-1/M^2)\|\Delta P_n x - P_n V'(P_n x)\|^2}{1+(1/M^2)\mathcal{E}_n(x)} \\ &\quad + \frac{1}{2}\left(1-\frac{1}{n}\right)\left\|\frac{(1/M^2)B_n^*(x)(\Delta P_n x - P_n V'(P_n x))}{1+(1/M^2)\mathcal{E}_n(x)}\right\|^2 \\ &\quad + \frac{1}{2n}\frac{\sum_{k=(1,\dots,1)}^{(m_n,\dots,m_n)}\sum_{i\equiv(i_1,\dots,i_d)=(1,\dots,1)}^{(m_n,\dots,m_n)}\lambda_i(\langle (B_n(x)/M^2)e_k, e_i\rangle)^2}{1+(1/M^2)\mathcal{E}_n(x)} \\ &\quad + \frac{1}{2n}\left(\sum_{k=(1,\dots,1)}^{(m_n,\dots,m_n)}\sum_{i=(1,\dots,1)}^{(m_n,\dots,m_n)}\sum_{j=(1,\dots,1)}^{(m_n,\dots,m_n)}\left\langle\frac{\sigma(\cdot,x)}{M^2}e_k, e_j\right\rangle\langle V''(P_n x)e_j, e_i\rangle \right. \\ &\qquad\qquad \left. \times \left\langle\frac{\sigma(\cdot,x)}{M^2}e_k, e_i\right\rangle\right) \end{aligned} \tag{A.9}$$

$$\begin{aligned}
&\quad \times \left(1 + \frac{1}{M^2}\mathcal{E}_n(x)\right)^{-1} \\
&\le \left\{ -\frac{1}{1 + (1/M^2)\mathcal{E}_n(x)} + \frac{1}{2}\left(1 - \frac{1}{n}\right)\left(\frac{1}{1 + (1/M^2)\mathcal{E}_n(x)}\right)^2 \right\} \\
&\quad \times \frac{\|\Delta P_n x - P_n V'(P_n x)\|^2}{M^2} \\
&\quad + \frac{1}{2n} m_n^d \sum_{i=(1,\dots,1)}^{(m_n,\dots,m_n)} \lambda_i + \frac{1}{n} m_n^{3d} \int_{\mathcal{O}} |V''(P_n x(\theta))| \, d\theta \\
&\le 0 + \frac{4^d \pi^{2d} (1 + m_n)^{4d}}{6n} + \frac{m_n^{3d}}{n} \sup_r |V''(r)| \le \text{Constant} < \infty,
\end{aligned}$$

where the $\lambda_i$'s are eigenvalues defined in (1.12) and the constant is independent of $n$. In the last inequality above, we used (1.14), and the estimate that

$$\sum_{i=(1,\dots,1)}^{(m_n,\dots,m_n)} \lambda_i = \left(\sum_{i_1=1}^{m_n} \mu_{i_1}\right)^d \le 4^d \pi^{2d} \frac{(1 + m_n)^{3d}}{3}, \tag{A.10}$$

where $\mu_i$ is the one in (1.10).

Let

$$\tau_n \equiv \inf\{t > 0 : \mathcal{E}_n(X_n(t)) \ge C_1\}.$$

By optional sampling theorem

$$\begin{aligned}
&\sup_{x : \mathcal{E}(x) \le C_0} P(\mathcal{E}_n(X_n(t)) > C_1, \text{ some } 0 < t \le T | X_n(0) = x) \\
&\qquad \times e^{n(C_1 - C_0) - nT \sup_{n,x} H_n f_n(x)} \\
&\quad \le \sup_{x : \mathcal{E}(x) \le C_0} E[e^{n f_n(X_n(T \wedge \tau_n)) - n f_n(X_n(0)) - \int_0^{T \wedge \tau_n} n H_n f_n(X_n(s)) \, ds} | X_n(0) = x] \\
&\quad = 1.
\end{aligned} \tag{A.11}$$

Hence (A.6) follows. We now relax the initial condition in the estimate to that in Lemma A.5/Condition 1.12. We achieve this by the following result:

LEMMA A.6. *We denote by $X_n^x$ the solution of (A.5) with initial value $X_n(0) = x$. Then for each $T, a > 0$, there exists a constant $C = C(T, a) > 0$ such that*

$$P\left(\sup_{0 \le t \le T} \|X_n^x(t) - X_n^y(t)\| > C\varepsilon \Big| \|x - y\| < \varepsilon\right) < e^{-na}$$

$$\forall\, 0 < \varepsilon < 1, n = 1, 2, \dots.$$

PROOF. For each $n$ fixed, $(X_n^x(t), X_n^y(t))$ is a two-component Markov process that solves the martingale problem with generator

$$\begin{aligned}\mathcal{A}_n f(x,y) = &\langle \Delta P_n x - P_n V'(P_n x), D_x f(x,y)\rangle \\ &+ \langle \Delta P_n y - P_n V'(P_n y), D_y f(x,y)\rangle \\ &+ \frac{1}{2n}\operatorname{Tr}((D_{xx}^2 f)B_n(x)B_n^*(x) \\ &+ (D_{yy}^2 f)B_n(y)B_n^*(y) + 2(D_{xy}^2 f)B_n(x)B_n^*(y))\end{aligned}$$

for $f(x,y) \in C^2(L^2(\mathcal{O}) \times L^2(\mathcal{O}))$. Let $\varepsilon > 0$ and

$$f_{n,\varepsilon}(x,y) \equiv \log\left(1 + \frac{1}{2}\left\|\frac{P_n x - P_n y}{\varepsilon}\right\|^2\right).$$

It follows that

$$\mathcal{H}_n f_{n,\varepsilon}(x,y) \equiv \frac{1}{n} e^{-n f_{n,\varepsilon}} \mathcal{A}_n e^{n f_{n,\varepsilon}}(x,y) \le C_0 < \infty,$$

where constant $C_0$ is independent of $n$ as well as $\varepsilon$. By an argument identical to that used in the proof of (A.11), the conclusion follows. □

PROOF OF LEMMA A.5. Let compact set $K \subset L^2(\mathcal{O})$ and $a, T, \varepsilon > 0$. It is enough for us to show that for any $x_n \in K$, there exists compact set $K_1 \subset E$,

$$P(\exists t \in [0,T], X_n^{x_n}(t) \notin K_1^{2\varepsilon}) \le 2e^{-na}.$$

Let $\delta > 0$. By compactness of $K$, there exists $\{x_0^1, \dots, x_0^{m(\delta)}\}$ such that

$$K \subset \bigcup_{k=1}^{m(\delta)} B(x_0^k, \delta).$$

Since $\{x : \mathcal{E}(x) < +\infty\} \subset H^1(\mathcal{O})$ is dense in $L^2(\mathcal{O})$, we can select $x_0^k$ so that

$$\sup_{j=1,\dots,m(\delta)} \mathcal{E}(x_0^j) < +\infty.$$

Therefore, we can choose

$$x_{0,n} \in \{x_0^{(1)}, \dots, x_0^{(m(\delta))}\}$$

such that $\|x_{0,n} - x_n\| < \delta$. By Lemma A.6, there exists $C = C(T,a) > 0$ (independent of $\delta$) such that

$$P\left(\sup_{0\le t\le T} \|X_n^{x_{0,n}}(t) - X_n^{x_n}(t)\| > C\delta\right) < e^{-na}.$$

By (A.6) and the compactness of level sets for $\mathcal{E}$, there exists a compact set $K_1 \subset E$ such that

$$P(X_n^{x_{0,n}}(t) \notin K_1^\delta, \exists\, t \in [0,T]) \le e^{-na}.$$

It follows that

$$\begin{aligned}&\{\exists\, t \in [0,T], X_n^{x_n}(t) \notin K_1^{(1+C)\delta}\}\\&\qquad \subset \{\exists\, t \in [0,T], X_n^{x_{0,n}}(t) \notin K_1^\delta\} \cup \left\{\sup_{0\le t\le T} \|X_n^{x_{n,0}}(t) - X_n^{x_n}(t)\| > C\delta\right\}.\end{aligned}$$

Therefore

$$\begin{aligned}&P(\exists\, t \in [0,T], X_n^{x_n}(t) \notin K_1^{(1+C)\delta})\\&\qquad \le P(\exists\, t \in [0,T], X_n^{x_{0,n}}(t) \notin K_1^\delta) + P\left(\sup_{0\le t\le T} \|X_n^{x_{0,n}}(t) - X_n^{x_n}(t)\| > C\delta\right)\\&\qquad \le 2e^{-na}.\end{aligned}$$

Taking $\delta = \varepsilon/(1+C)$, we complete the proof. $\square$

A.2.2. *Stochastic Cahn–Hilliard equation.* Recall the stochastic Cahn–Hilliard equation (1.22) in Example 1.5:

$$dX_n(t) = \Delta P_n(-\Delta P_n X_n(t) + P_n V'(P_n X_n(t)))\, dt + \frac{1}{\sqrt{n}} B_n\, dW(t). \tag{A.12}$$

We prove the following lemma.

LEMMA A.7. *Condition 1.12 holds for $X_n$.*

Using identical arguments as in the Allen–Cahn case, we just need the following estimates: $\sup_n \sup_x H_n f_n(x) < \infty$ for the $f_n$ below, and (A.15).

Define $\mathcal{E}, \mathcal{E}_n$ the same way as in (A.7). Let

$$f_n(x) \equiv \log\left(1 + \frac{1}{M^2}\mathcal{E}_n(x)\right),$$

where $M > 0$ is the constant in the Poincaré type inequality

$$\left\| x - \int_{\theta\in\mathcal{O}} x(\theta)\, d\theta \right\| \le M\|\nabla x\| \qquad \forall\, x \in H^2(\mathcal{O}). \tag{A.13}$$

Let $\lambda_k$ be defined according to (1.12). Then

$$\begin{aligned}H_n f_n(x) &\equiv \langle \Delta P_n(-\Delta P_n x + P_n V'(P_n x)), Df_n(x)\rangle\\&\qquad + \frac{1}{2}\|Df_n(x)\|^2 + \frac{1}{2n}\operatorname{Tr}(D^2 f_n(x))\end{aligned}$$

$$
\begin{aligned}
&= \frac{(-1/M^2)\|(-\Delta P_n)^{1/2}(\Delta P_n x - P_n V'(P_n x))\|^2}{1+(1/M^2)\mathcal{E}_n(x)} \\
&\quad + \frac{1}{2}\left(1-\frac{1}{n}\right)\frac{\|\Delta P_n x - P_n V'(P_n x)\|^2/M^4}{(1+(1/M^2)\mathcal{E}_n(x))^2} \\
&\quad + \frac{1}{2n}\frac{(1/M^2)\sum_{k=(1,\dots,1)}^{(m_n,\dots,m_n)}\langle(-\Delta P_n + P_n V''(P_n x))e_k, e_k\rangle}{1+(1/M^2)\mathcal{E}_n(x)} \\
&= \frac{(-1/M^2)\|\nabla(\Delta P_n x - P_n V'(P_n x))\|^2}{1+(1/M^2)\mathcal{E}_n(x)} \\
&\quad + \frac{1}{2}\left(1-\frac{1}{n}\right)\frac{(1/M^2)\|\Delta P_n x - P_n V'(P_n x)\|^2/M^2}{(1+(1/M^2)\mathcal{E}_n(x))^2} \\
&\quad + \frac{1}{M^2}\frac{1}{2n}\frac{\sum_{k=(1,\dots,1)}^{(m_n,\dots,m_n)}\lambda_k + \sum_{k=(1,\dots,1)}^{(m_n,\dots,m_n)}\langle V''(P_n x)e_k, e_k\rangle}{1+(1/M^2)\mathcal{E}_n(x)} \\
&\le \left\{-\frac{1}{1+(1/M^2)\mathcal{E}_n(x)} + \frac{1}{2}\left(1-\frac{1}{n}\right)\left(\frac{1}{1+(1/M^2)\mathcal{E}_n(x)}\right)^2\right\} \qquad \text{(A.14)} \\
&\quad \times \frac{\|\nabla(\Delta P_n x - P_n V'(P_n x))\|^2}{M^2} \\
&\quad + \frac{1}{2}\left(\frac{1}{M^2}\int_{\mathcal{O}} V'(P_n x(\theta))\,d\theta\right)^2 \\
&\quad + \frac{1}{2n}\frac{1}{M^2}\sum_{k=(1,\dots,1)}^{(m_n,\dots,m_n)}\left(\lambda_k + \sup_r |V''(r)|\right) \\
&\le C < \infty.
\end{aligned}
$$

In the above derivations, we used (A.13):

$$
\begin{aligned}
&\|\Delta P_n x - P_n V'(P_n x)\|^2 \\
&\quad = \left\{\int_{\mathcal{O}}(\Delta P_n x(\theta) - P_n V'(P_n x(\theta)))\,d\theta\right\}^2 \\
&\qquad + \left\|\Delta P_n x - P_n V'(P_n x) - \int_{\mathcal{O}}(\Delta P_n x(\theta) - P_n V'(P_n x(\theta)))\,d\theta\right\|^2 \\
&\quad \le \left(\int_{\mathcal{O}} V'(P_n x(\theta))\,d\theta\right)^2 + M^2\|\nabla(\Delta P_n x - P_n V'(P_n x))\|^2.
\end{aligned}
$$

We also made use of (A.10) and condition (1.21).

LEMMA A.8. *We denote by $X_n^x$ the solution of (A.12) with initial value $X_n(0) = x$. Then for each $T, a > 0$, there exists a constant $C = C(T, a) > 0$*

*such that*

$$P\Big(\sup_{0\le t\le T}\|X_n^x(t)-X_n^y(t)\|>C\varepsilon\Big|\|x-y\|<\varepsilon\Big)<e^{-na} \qquad \forall\, 0<\varepsilon<1, n=1,2,\ldots. \tag{A.15}$$

Proof. The proof follows the same idea as in Lemma A.6. □

A.2.3. *Stochastic quasilinear equation with viscosity.* Using the same ideas as in the stochastic Allen–Cahn and Cahn–Hilliard case, by choosing

$$\mathcal{E}(x)=\tfrac{1}{2}(\|x\|^2+\|\nabla x\|^2),\qquad \mathcal{E}_n(x)=\tfrac{1}{2}(\|P_nx\|^2-\langle\Delta P_nx,x\rangle)$$

and

$$f_n(x)=\log(1+\alpha\mathcal{E}_n(x)),$$

we can prove the following.

Lemma A.9. *Condition 1.12 holds for the $X_n$ in (1.26).*

**A.3. Approximations of the Tataru distance function.** Let $E,U_0$ be real separable Hilbert spaces. We discuss approximations of the *Tataru distance function* $d_C$ (Definition 1.13) by $C^2(E)$ functions. Throughout this section, we assume Condition 1.11 is satisfied for $C,C_n$.

We want to keep two useful properties (3.3) and (3.4) in the approximation. The functions $h_{\varepsilon,y}(x)$ and $h_{n,\varepsilon,y}(x)$ defined in (3.8) and (3.9) satisfy these requirements—see (A.18) and (A.20).

Lemma A.10. *For each $\varepsilon>0$ small enough, define $\phi_\varepsilon$ according to (3.7): $\phi_\varepsilon(r)=\sqrt{r}$ when $r\ge\varepsilon$ and*

$$\phi_\varepsilon(r)=\sqrt{\varepsilon}+\frac{r-\varepsilon}{2\sqrt{\varepsilon}}-\frac{(r-\varepsilon)^2}{8\varepsilon\sqrt{\varepsilon}}\qquad \text{when } 0\le r\le\varepsilon.$$

*Then:*

(1) $\phi_\varepsilon',\phi_\varepsilon''\in C_b([0,+\infty))$; *$\phi_\varepsilon$ is nondecreasing,* $\sup_r r|\phi_\varepsilon''(r)|<+\infty$.
(2) $\lim_{\varepsilon\to0}\sup_{r\ge0}|\phi_\varepsilon(r)-\sqrt{r}|=0$.
(3)

$$0\le r\phi_\varepsilon'(r^2)\le 1/2. \tag{A.16}$$

Proof. The first two properties follow from direct verification. To see that the third one holds, let $f(r)=r\phi_\varepsilon'(r^2)$; then $f'(r)>0$ for $0<r<\varepsilon$, and $f'(r)=0$ when $r\ge\varepsilon$. Hence $r\phi_\varepsilon'(r^2)\le\varepsilon\phi_\varepsilon'(\varepsilon^2)=1/2$. □

Lemma A.11. *Let $\varphi, \varphi_n : [0,\infty) \to [0,\infty)$ be continuous. Suppose there exist $0 < m < M < +\infty$, $0 \le c < \infty$ such that $mt \le \varphi_n(t) \le c + Mt$, $n = 1, 2, \ldots$. Suppose further that $0 < a_n \to +\infty$ and that*

$$\lim_{n\to\infty} \sup_{0\le t\le T} |\varphi_n(t) - \varphi(t)| = 0 \qquad \forall T \ge 0.$$

*Then*

$$\lim_{n\to\infty} -\frac{1}{a_n} \log \int_0^\infty e^{-a_n \varphi_n(t)}\, dt = \inf_{t\ge 0} \varphi(t).$$

Proof. It is straightforward to see that, when $T > 0$ is large enough but fixed,

$$\lim_{n\to\infty} -\frac{1}{a_n} \log \int_0^T e^{-a_n\varphi_n(t)}\, dt = \inf_{0\le t\le T} \varphi(t) = \inf_{t\ge 0} \varphi(t).$$

Take $T \ge \varphi(0)/m$; then when $n$ is large enough

$$-\frac{1}{a_n} \log \int_T^\infty e^{-a_n\varphi_n(t)}\, dt \ge -\frac{1}{a_n} \log \int_T^\infty e^{-a_n m t}\, dt \ge mT \ge \varphi(0) \ge \inf_{t\ge 0} \varphi(t).$$

Therefore,

$$\begin{aligned}
\inf_{t\ge 0} \varphi(t) &= \min\left\{ \liminf_{n\to\infty} -\frac{1}{a_n} \log \int_0^T e^{-a_n\varphi_n(t)}\, dt, \liminf_{n\to\infty} -\frac{1}{a_n} \log \int_T^\infty e^{-a_n\varphi_n(t)}\, dt \right\} \\
&\le \liminf_{n\to\infty} -\frac{1}{a_n} \log \int_0^\infty e^{-a_n\varphi_n(t)}\, dt \le \limsup_{n\to\infty} -\frac{1}{a_n} \log \int_0^\infty e^{-a_n\varphi_n(t)}\, dt \\
&\le \limsup_{n\to\infty} -\frac{1}{a_n} \log \int_0^T e^{-a_n\varphi_n(t)}\, dt = \inf_{t\ge 0} \varphi(t).
\end{aligned}$$

□

Lemma A.12. *Let $a_n > 0$ be such that $\lim_{n\to\infty} a_n = \infty$. Define $\phi_\varepsilon$ as in (3.7) and $h_\varepsilon, h_{n,\varepsilon}$ according to (3.8) and (3.9):*

$$h_{n,\varepsilon}(x) \equiv -\frac{1}{a_n} \log \int_0^\infty e^{-a_n\{t+\phi_\varepsilon(\|x-S_n(t)y\|^2)\}}\, dt$$

*and*

$$h_\varepsilon(x) \equiv \inf_{t\ge 0} \{t + \phi_\varepsilon(\|x - S(t)y\|^2)\}.$$

*Then:*

(1) *$h_{n,\varepsilon}(x) \ge c$ whenever $\|x\| \ge \|y\| + c$, for every $c > 1$, $a_n > 1$.*
(2) *For each $y \in E$ fixed,*

$$\lim_{\varepsilon\to 0+} \sup_{x\in E} |h_\varepsilon(x) - d_C(x,y)| = 0. \tag{A.17}$$

(3) *For each* $\varepsilon > 0$ *fixed,*

$$\lim_{n\to\infty} \sup_{x\in K} |h_\varepsilon(x) - h_{n,\varepsilon}(x)| = 0 \qquad \forall \text{ compact } K \subset E.$$

(4) $h_{n,\varepsilon} \in C^2(E)$;

$$\|Dh_{n,\varepsilon}(x)\| \le 1. \tag{A.18}$$

*If* $B_n$ *is a Hilbert–Schmidt operator from* $U_0$ *to* $E$ [*i.e.,* $B_n \in L_2(U_0, E)$], *then*

$$\begin{aligned} &|\operatorname{Tr}[D^2 h_{n,\varepsilon}(x) B_n B_n^*]| \\ &\quad \le \Big(2a_n + 4\sup_{r\ge 0} r|\phi_\varepsilon''(r)| + 2\sup_{r\ge 0}\phi_\varepsilon'(r)\Big) \|\|B_n\|\|^2_{L_2(U_0,E)}. \end{aligned} \tag{A.19}$$

(5) *Let* $\varphi \in \mathcal{T}$; *we have*

$$\limsup_{r\to 0+, z\to x} \frac{\varphi(h_{n,\varepsilon}(S_n(r)z)) - \varphi(h_{n,\varepsilon}(z))}{r} \le \sup_{r\ge 0}\varphi'(r) \qquad \forall x \in E. \tag{A.20}$$

*If* $B_n \in L_2(U_0, E)$, *then*

$$\begin{aligned} &|\operatorname{Tr}[D^2(\varphi\circ h_{n,\varepsilon})(x) B_n B_n^*]| \\ &\quad \le \Big\{\sup_{s\ge 0}\varphi''(s) + \sup_{s\ge 0}\varphi'(s)\Big(2a_n + 4\sup_{r\ge 0} r|\phi_\varepsilon''(r)| + 2\sup_{r\ge 0}\phi_\varepsilon'(r)\Big)\Big\} \\ &\quad\quad \times \|\|B_n\|\|^2_{L_2(U_0,E)}. \end{aligned} \tag{A.21}$$

PROOF. Part (1): Since $\|S_n(t)y\| \le \|y\|$, $(\|x\| - \|y\|) \vee 0 \le \|x - S_n(t)y\|$. Hence

$$\phi_\varepsilon(((\|x\| - \|y\|)\vee 0)^2) \le h_{n,\varepsilon}(x)$$

when $a_n > 1$. Noting $\phi_\varepsilon(r^2) = r$ when $r \ge 1$, the conclusion follows.

Part (2) is a direct consequence of part (2) of Lemma A.10.

Part (3) follows if we prove that for each $x_n \to x$,

$$\lim_{n\to\infty} h_{n,\varepsilon}(x_n) = h_\varepsilon(x).$$

Take $\varphi_n(t) = \{t + \phi_\varepsilon(\|x_n - S_n(t)y\|^2)\}$, $\varphi(t) = \{t + \phi_\varepsilon(\|x - S(t)y\|^2)\}$. For each $T > 0$,

$$\begin{aligned} \sup_{0\le t\le T} |\varphi_n(t) - \varphi(t)| &= \sup_{0\le t\le T} |\phi_\varepsilon(\|x_n - S_n(t)y\|^2) - \phi_\varepsilon(\|x - S(t)y\|^2)| \\ &\le \sup_{r\ge 0}\phi_\varepsilon'(r) \sup_{0\le t\le T} |\|x_n - S_n(t)y\| - \|x - S(t)y\|| \to 0 \\ &\qquad \text{as } n\to\infty. \end{aligned}$$

In addition

$$t \le \varphi_n(t) \le t + \sup_n \phi_\varepsilon((\|x_n\| + \|y\|)^2).$$

Apply Lemma A.11; therefore

$$\lim_{n\to\infty} -\frac{1}{a_n} \log \int_0^\infty e^{-a_n \varphi_n(t)}\, dt = \inf_{t\ge 0} \varphi(t).$$

Part (4): It can be verified that

$$Dh_{n,\varepsilon}(x) = 2 \frac{\int_0^\infty e^{-a_n\{t+\phi_\varepsilon(\|x-S_n(t)y\|^2)\}} \phi'_\varepsilon(\|x-S_n(t)y\|^2)(x-S_n(t)y)\, dt}{\int_0^\infty e^{-a_n\{t+\phi_\varepsilon(\|x-S_n(t)y\|^2)\}}\, dt} \in E.$$

By Lemma A.10,

$$\|Dh_{n,\varepsilon}(x)\| \le 1.$$

Direct calculation also gives

$$\begin{aligned} D^2 h_{n,\varepsilon}(x) &= \frac{1}{\int_0^\infty e^{-a_n\{t+\phi_\varepsilon(\|x-S_n(t)y\|^2)\}}\, dt} \\ &\quad \times \int_0^\infty e^{-a_n\{t+\phi_\varepsilon(\|x-S_n(t)y\|^2)\}} \\ &\qquad\qquad \times ((-4a_n)(\phi'_\varepsilon(\|x-S_n(t)y\|^2))^2 \\ &\qquad\qquad\quad \times (x-S_n(t)y)\otimes(x-S_n(t)y) \\ &\qquad\qquad\quad + 4\phi''_\varepsilon(\|x-S_n(t)y\|^2)(x-S_n(t)y)\otimes(x-S_n(t)y) \\ &\qquad\qquad\qquad\qquad\qquad + 2\phi'_\varepsilon(\|x-S_n(t)y\|^2) I)\, dt \\ &\quad + a_n Dh_{n,\varepsilon}(x)\otimes Dh_{n,\varepsilon}(x). \end{aligned}$$

Let $\{\hat{e}_1, \ldots, \hat{e}_k, \ldots\}$ be a complete orthonormal basis for $E$. Then,

$$\begin{aligned} &|\operatorname{Tr}[D^2 h_{n,\varepsilon}(x) B_n B_n^*]| \\ &\quad = \left|\sum_k \langle D^2 h_{n,\varepsilon}(x) B_n B_n^* \hat{e}_k, \hat{e}_k\rangle\right| \\ &\quad \le \sup_{t\ge 0}(4a_n(\phi'_\varepsilon(\|x-S_n(t)y\|^2))^2 + 4|\phi''_\varepsilon|(\|x-S_n(t)y\|^2)) \\ &\qquad\quad \times \|B_n^*(x-S_n(t)y)\|_{U_0}^2 \\ &\qquad + 2\sup_{t\ge 0} \phi'_\varepsilon(\|x-S_n(t)y\|^2)\|B_n\|_{L_2(U_0,E)}^2 + a_n\|B_n^* Dh_{n,\varepsilon}(x)\|^2 \\ &\quad \le 4a_n \sup_{t\ge 0}(\|x-S_n(t)y\|\phi'_\varepsilon(\|x-S_n(t)y\|^2)\|B_n\|_{L_2(U_0,E)})^2 \end{aligned}$$

$$
\begin{aligned}
&+ 4 \sup_{t\geq 0} |\phi''_\varepsilon|(\|x - S_n(t)y\|^2)\|x - S_n(t)y\|^2 \|\!|B_n|\!\|^2_{L_2(U_0,E)} \\
&+ 2 \sup_{t\geq 0} \phi'_\varepsilon(\|x - S_n(t)y\|^2)\|\!|B_n|\!\|^2_{L_2(U_0,E)} + a_n \|Dh_{n,\varepsilon}(x)\|^2 \|\!|B_n|\!\|^2_{L_2(U_0,E)} \\
\leq{}& \Big(4a_n \sup_{r\geq 0}(r\phi'_\varepsilon(r^2))^2 + 4 \sup_{r\geq 0} r|\phi''_\varepsilon(r)| + 2 \sup_{r\geq 0} \phi'_\varepsilon(r) + a_n\Big)\|\!|B_n|\!\|^2_{L_2(U_0,E)}.
\end{aligned}
$$

To derive the second inequality above, we used the following: let $e_1, \dots, e_k, \dots$ be a complete orthonormal system for $U_0$; then for each $z \in E$,

$$
\begin{aligned}
\|B_n^* z\|^2 &= \sum_k \langle B_n^* z, e_k\rangle^2 = \sum_k \langle z, B_n e_k\rangle^2 \leq \sum_k \|z\|^2 \|B_n e_k\|^2 \\
&= \|z\|^2 \|\!|B_n|\!\|^2_{L_2(U_0,E)}.
\end{aligned}
$$

Noting (A.16), (A.19) holds.

Part (5): Let $r > 0$; then

$$
\begin{aligned}
h_{n,\varepsilon}(x) + r &= -\frac{1}{a_n} \log \int_0^\infty e^{-a_n(t+r+\phi_\varepsilon(\|x - S_n(t)y\|^2))}\,dt \\
&\geq -\frac{1}{a_n} \log \int_0^\infty e^{-a_n(t+r+\phi_\varepsilon(\|S_n(r)x - S_n(t+r)y\|^2))}\,dt \\
&= -\frac{1}{a_n} \log \int_r^\infty e^{-a_n(t+\phi_\varepsilon(\|S_n(r)x - S_n(t)y\|^2))}\,dt \\
&\geq h_{n,\varepsilon}(S_n(r)x).
\end{aligned}
$$

Hence we have (A.20). Equation (A.21) follows from (A.18), (A.19) and

$$
D^2\varphi \circ h_{n,\varepsilon}(x) = \varphi'' \circ h_{n,\varepsilon}(x)(Dh_{n,\varepsilon}(x) \otimes Dh_{n,\varepsilon}(x)) + \varphi' \circ h_{n,\varepsilon}(x) D^2 h_{n,\varepsilon}(x).
$$

□

**Acknowledgments.** I would like to thank Professor Tom Kurtz for useful discussions on stochastic equation in infinite dimensions, and Professor Markos Katsoulakis for useful discussions on stochastic Allen–Cahn and Cahn–Hilliard equations. An anonymous referee carefully read through the paper, and provided an extensive list of thoughtful comments and suggestions. These are incorporated in the paper; I gratefully acknowledge the help as well.

Department of Mathematics and Statistics
University of Massachusetts–Amherst
Amherst, Massachusetts 01003
USA
e-mail: feng@math.umass.edu